\documentclass[11pt]{article}
\usepackage{amssymb,amsmath,euscript}
\newtheorem{proposition}{Proposition}[section]

\newtheorem{lemma}[proposition]{Lemma}
\newtheorem{theorem}[proposition]{Theorem}

\newtheorem{corollary}[proposition]{Corollary}

\def\l{{\langle}}
\def\r{\rangle}
\def\dim{{\rm dim}_{_{\rm H}}}

\def\dimp{{\rm dim}_{_{\rm P}}}

\def\R{{\mathbb R}}

\def\a{\alpha}
\def\la{\lambda}

\def\ga{\gamma}
\def\Ga{\Gamma}
\def\ep{\varepsilon}
\def\eps{\varepsilon}
\def\si{\sigma}

\def\t{{\bf t}}

\def\E{{\mathbb E}}
\def\P{{\mathbb P}}

\makeatletter \@addtoreset{equation}{section} \makeatother
\newcommand{\qed}%
{%
    {}\hfill
    {}\hfill
    {$\square $}%
    \vspace {0.3cm}%
    \pagebreak [2]%
    \par
}%
\newenvironment{proof}[1]{%
    \vspace{0.3cm}%
    \pagebreak [2]%
    \par%
    \noindent {\bf  Proof~#1\ }}{\qed}%
\newenvironment{remark}{%
    \vspace{0.3cm} \pagebreak [2]%
    \par%
    \refstepcounter{proposition}
    \noindent%
   {\bf Remark~\theproposition\  }}{\ }%
\begin{document}
\date{}

\title {Regularity of Intersection Local Times of Fractional Brownian Motions}
\author{Dongsheng Wu\\ University of Alabama in Huntsville\\ \and
Yimin Xiao \footnote{Research partially
        supported by  NSF grant DMS-0706728.}  \\
Michigan State University }

\maketitle

\begin{abstract}
Let $B^{\alpha_i}$ be an $(N_i,d)$-fractional Brownian motion with
Hurst index ${\alpha_i}$ ($i=1,2$), and let $B^{\alpha_1}$ and
$B^{\alpha_2}$ be independent. We prove that, if
$\frac{N_1}{\alpha_1}+\frac{N_2}{\alpha_2}>d$, then the intersection
local times of $B^{\alpha_1}$ and $B^{\alpha_2}$ exist, and have a
continuous version. We also establish H\"{o}lder conditions for the
intersection local times and determine the Hausdorff and packing
dimensions of the sets of intersection times and intersection
points.

One of the main motivations of this paper is from the results of
Nualart and Ortiz-Latorre ({\it J. Theor. Probab.} {\bf 20} (2007)),
where the existence of the intersection local times of two
independent $(1,d)$-fractional Brownian motions with the same Hurst
index was studied by using a different method. Our results show that
anisotropy brings subtle differences into the analytic properties of
the intersection local times as well as rich geometric structures
into the sets of intersection times and intersection points.
\end{abstract}

      {Running head}:  Regularity of intersection local times of
      fractional Brownian motions\\

      {\it 2000 AMS Classification numbers}: 60G15, 60J55, 60G18, 60F25, 28A80.\\



\section{Introduction}\label{Sec:Intro}

Let $B_0^\gamma=\{B_0^\gamma(u),\,u\in\R^p\}$ be a $p$-parameter
fractional Brownian motion in $\R$ with Hurst index
$\gamma\in(0,\,1),$ i.e., a centered, real-valued Gaussian random
field with covariance function
\begin{equation}\label{Eq:Cov_fBm}
\E
\left[B_0^\gamma(u_1)B_0^\gamma(u_2)\right]=\frac12\left(|u_1|^{2\gamma}
+|u_2|^{2\gamma}-|u_1-u_2|^{2\gamma}\right).
\end{equation}
It follows from Eq. (\ref{Eq:Cov_fBm}) that
$\E\big[\big(B_0^\gamma(u_1)-B_0^\gamma(u_2)\big)^2\big]=|u_1-u_2|^{2\gamma}$
and $B_0^\gamma$ is $\gamma$-self-similar with stationary
increments.

We associate with $B_0^\gamma$ a Gaussian random field $B^\ga = \{B^\ga(u), u\in\R^p\}$ in
$\R^q$ by
\begin{equation}\label{Eq:pq_field}
B^\ga(u)=\big(B_1^\ga(u),\ldots,B_q^\ga(u)\big),\quad\, u\in\R^p,
\end{equation}
where $B_1^\ga,\ldots,B_q^\ga$ are independent copies of $B_0^\ga.$ $B^\ga$ is called
a $(p, q)$-fractional Brownian motion of index $\ga$.

Fractional Brownian motion has been intensively studied in recent
years and, because of its interesting properties such as short/long
range dependence and self-similarity, has been widely applied in
many areas such as finance, hydrology and telecommunication
engineering.

Let $B^{\a_1}=\{B^{\a_1}(s),\,s\in\R^{N_1}\}$ and
$B^{\a_2}=\{B^{\a_2}(t),\,t\in\R^{N_2}\}$ be two independent
fractional Brownian motions in $\R^d$ with Hurst indices $\a_1,\,
\a_2\in(0,\,1)$, respectively. This paper is concerned with the
regularity of the intersection local times of $B^{\a_1}$ and
$B^{\a_2}$, as well as the fractal properties of the sets of intersection
times and intersection points. Without loss of generality, we further
assume $\a_1\le \a_2$ throughout this paper. For $N_1=N_2=1$ and
${\a_1}={\a_2}=\frac12$, the processes are classical $d$-dimensional
Brownian motions. The intersection local times of independent
Brownian motions have been studied by several authors [see Wolpert
(1978a), Geman, Horowitz and Rosen (1984)] and is closely related to
the self-intersections (or multiple points) of Brownian motion. The
approach of these papers relies on the fact that the intersection
local times of independent Brownian motions can be seen as the local
times at zero of some Gaussian random field. For the applications of
the intersection local time theory for Brownian motions, we refer to
Wolpert (1978b) and LeGall (1985), among others.

The self-intersection local times of fractional Brownian motion were
studied by Rosen (1987) for the planar case, and by Hu and Nualart
(2005) for the multidimensional case. Very recently, Nualart and
Ortiz-Latorre (2007)  proved an existence result for the
intersection local times of two independent $d$-dimensional
fractional Brownian motions with the same Hurst index.

The aim of this paper is to show that the existence of the
intersection local times for two independent fractional Brownian
motions $B^{\a_1}$ and $B^{\a_2}$ in $\R^d$ can be studied by using
a Fourier analytic method and, moreover, this latter method can be
applied to establish the joint continuity and sharp H\"{o}lder
conditions for the intersection local times. Besides their own
interest, these results are useful for studying fractal properties
of the set of intersection times as well as the set of intersection
points.

Let $X=\{X(s,t),\,(s,t)\in\R^N\}$ be an $(N,d)$-Gaussian random
field, where $N=N_1+N_2$, defined by
\begin{equation}\label{Def:X}
X(s,t)\equiv B^{\a_1}(s)-B^{\a_2}(t),\qquad
s\in\R^{N_1},\,\,t\in\R^{N_2}.
\end{equation}
We will follow the same idea as Wolpert (1978a) and Geman, Horowitz
and Rosen (1984) and treat the intersection local times of
$B^{\a_1}$ and $B^{\a_2}$ as the local times at 0 of $X$, with an
intension to establish sharp H\"older conditions. The main
ingredients for proving our results are the strong local
nondeterminism of fractional Brownian motions, occupation density
theory [cf. Geman and Horowitz (1980)], and newly developed
techniques for anisotropic Gaussian random fields [cf. Ayache, Wu
and Xiao (2008) and Xiao (2009)].

For later use, we mention that, by the self-similarity and
stationarity of the increments of $B^{\a_1}$ and $B^{\a_2}$,
the Gaussian random field $X$ defined by (\ref{Def:X}) has stationary
increments and satisfies the following
operator-scaling property: For every constant $c > 0$,
\begin{equation}\label{Eq:OSS}
\left\{X(c^A\,(s,t)), (s,t) \in \R^N\right\} \stackrel{d}{=} \left\{c\,X(s,t),
(s, t) \in \R^N\right\},
\end{equation}
where $A=(a_{ij})$ is an $N\times N$ diagonal matrix such that
$a_{ii} = 1/\a_1$ if $1 \le i \le N_1$ and $a_{ii} = 1/\a_2$
if $N_1+1 \le i \le N$. In the above, $\stackrel{d}{=}$ denotes equality
of all finite dimensional distributions and $c^A$ is the linear operator on $\R^N$
defined by $c^A = \sum_{n=0}^\infty \frac{(\ln c)^n A^n} {n!}.$

This paper is organized as follows. In Section \ref{Sec:Pre}, we
give several lemmas which will be used to prove our main results in
the following sections. In Section \ref{Sec:LT}, we study the
existence and the joint continuity of the intersection local times
of two independent $d$-dimensional fractional Brownian motions. We
prove that the necessary and sufficient condition for the existence
of an intersection local times in $L^2(\P\times \lambda_d)$ actually
implies the joint continuity. We devote Section \ref{Sec:H} to the
study of the exponential integrability and H\"{o}lder conditions for the
intersection local times. The later results imply information about
the exact Hausdorff measure of the set of intersection times of
$B^{\a_1}$ and $B^{\a_2}$. Finally, in Section \ref{Sec:dim}, we
determine the Hausdorff and packing dimensions of the set of intersection
points of $B^{\a_1}$ and $B^{\a_2}$.

Throughout this paper, we use $\l \cdot, \cdot\r$ and $|\cdot|$ to
denote the ordinary scalar product and the Euclidean norm in $\R^p$,
respectively, no matter what the value of the integer $p$ is, and we
use $\lambda_p$ to denote the Lebesgue measure in $\R^p$. We denote by
$O_p(u,r)$ a $p$-dimensional ball centered at $u$ with radius $r$,
and $O_{p_1,p_2}\big(u,r\big):=O_{p_1}(u_1,r)\times O_{p_2}(u_2,r),$
where $u=(u_1,u_2)$ with $u_1\in\R^{p_1}$ and $u_2\in\R^{p_2}$. In
Section $i$, unspecified positive and finite constants will be
numbered as $c_{_{i,1}},\, c_{_{i,2}},$....

\section{Preliminaries}\label{Sec:Pre}

In this section, we provide necessary preparations for the proofs of
our main results in the later sections.

It follows from Lemma 7.1 of Pitt (1978) that, for any $\gamma \in
(0, 1)$, the real-valued fractional Brownian motion $B_0^\ga =
\{B_0^\ga(u), u \in \R^p\}$ has the following important property of
\emph{strong local nondeterminism}: There exists a constant
$0<c_{_{2,1}}<\infty$ such that for all integers $n \ge 1$ and all
$u,\, u_1,\ldots,u_n \in \R^p$,
\begin{equation} \label{Eq:LND}
\mathrm{Var}\Big(B_0^\ga(u)\big|B_0^\ga(u_1),\ldots,
B_0^\ga(u_n)\Big) \ge c_{_{2,1}}\min_{0 \le k \le n}|u-u_k|^{2\ga},
\end{equation}
where $\mathrm{Var}\big(B_0^\ga(u)\big|B_0^\ga(u_1),\ldots,
B_0^\ga(u_n)\big)$ denotes the conditional variance of $B_0^\ga(u)$
given $B_0^\ga(u_1),$ $\ldots, B_0^\ga(u_n)$, and where $u_0 \equiv
0$.  The strong local nondeterminism has played important r\^o{}les
in studying various sample path properties of fractional Brownian
motion. See Xiao (1997, 2006, 2007) and the references therein for further
information. It will be the main technical tool of this paper as
well.

We consider the real-valued Gaussian random field $X_0 = \{X_0(s,t),
(s, t)\in \R^{N}\} $ defined by
$X_0(s,t):=B_0^{\a_1}(s)-B_0^{\a_2}(t)$ for $s\in \R^{N_1}$ and $t
\in \R^{N_2}$. Then the coordinate processes of $X$ defined by (\ref{Def:X})
are independent copies of $X_0$.

The following Lemma \ref{Lem:XLND} is a consequence
of the property of strong local nondeterminism of fractional
Brownian motion, and will be useful in our approach.

\begin{lemma}\label{Lem:XLND}
There exists a constant $0<c_{_{2,2}}<\infty$ such that
for all integers $n \ge 1$ and all
$(v,w),\,(s_1,t_1),\ldots,(s_n,t_n)\in\R^N$, we have
\begin{equation}\label{Eq:XLND}
\begin{split}
\mathrm{Var}&\Big(X_0(v,w)\big|X_0(s_1,t_1),\ldots,X_0(s_n,t_n)\Big)
\ge c_{_{2,2}}\, \bigg( \min_{0 \le k \le n}|v-s_k|^{2\a_1}+
\min_{0 \le k \le n}|w-t_k|^{2\a_2} \bigg),
\end{split}
\end{equation}
where $s_0 = t_0 = 0$.
\end{lemma}

\begin{proof}\ By definition we can write
\begin{equation}\label{Eq:XLND1}
\begin{split}
\mathrm{Var}\Big(X_0(v,w)\big|X_0(s_1,t_1),\ldots,X_0(s_n,t_n)\Big)
=\inf_{a_i\in\R,\,1 \le i \le n} \E\Bigg[\bigg(X_0(v,w)-
\sum_{i=1}^na_iX_0(s_i,t_i)\bigg)^2\Bigg].
\end{split}
\end{equation}
Since $B_0^{\a_1}$ and $B_0^{\a_2}$ are independent, we have
\begin{equation}\label{Eq:XLND1b}
\begin{split}
&\mathrm{Var}\Big(X_0(v,w)\big|X_0(s_1,t_1),\ldots,X_0(s_n,t_n)\Big)\\
&=\inf_{a_i\in\R,\, 1 \le i \le n}\left\{\E\Bigg[\bigg(B^{\a_1}_0(v)
-\sum_{i=1}^na_iB^{\a_1}_0(s_i)\bigg)^2\Bigg]
+\E\Bigg[\bigg(B^{\a_2}_0(w)-
\sum_{i=1}^na_iB^{\a_2}_0(t_i)\bigg)^2\Bigg]\right\}\\
&\ge \inf_{a_i\in\R,\, 1 \le i \le n}\E\Bigg[\bigg(B^{\a_1}_0(v)
-\sum_{i=1}^na_iB^{\a_1}_0(s_i)\bigg)^2\Bigg] +\inf_{b_i\in\R,\,1
\le i \le n}\E\Bigg[\bigg(B^{\a_2}_0(w)
-\sum_{i=1}^nb_iB^{\a_2}_0(t_i)\bigg)^2\Bigg]\\
&=\mathrm{Var}\Big(B^{\a_1}_0(v)\big|B^{\a_1}_0(s_1),\ldots,B^{\a_1}_0(s_n)\Big)+
\mathrm{Var}\Big(B^{\a_2}_0(w)\big|B^{\a_2}_0(t_1),\ldots,B^{\a_2}_0(t_n)\Big).
\end{split}
\end{equation}
Hence (\ref{Eq:XLND}) follows from (\ref{Eq:XLND1b}) and
(\ref{Eq:LND}).
\end{proof}

Combining Lemma \ref{Lem:XLND} with the following well-known fact,
which will be used repeatedly throughout the paper, that
\begin{equation}\label{Eq:Det}
{\rm detCov} (Z_1, \ldots, Z_n) = {\rm Var} (Z_1)\prod_{k=2}^n {\rm
Var}(Z_k | Z_1,\ldots, Z_{k-1})
\end{equation}
for any  Gaussian random vector $(Z_1, \ldots, Z_n)$, we have that,
for any $(s_1,t_1),\ldots,(s_n,t_n)\in\R^N_+,$
\begin{equation}\label{Eq:Det_XB}
\begin{split}
&{\rm detCov}\left(X_0(s_1,t_1),\ldots,X_0(s_n,t_n)\right)\\
&\ge  \prod_{j=1}^n \bigg[\mathrm{Var}
\big(B^{\a_1}_0(s_j)\big|B^{\a_1}_0(s_1),\ldots,B^{\a_1}_0(s_{j-1})\big)+
\mathrm{Var}\big(B^{\a_2}_0(t_j)\big|B^{\a_2}_0(t_1),\ldots,
B^{\a_2}_0(t_{j-1})\big)\bigg]\\
&\ge c_{_{2,2}}^n\, \prod_{j=1}^n  \bigg( \min_{0 \le k \le
j-1}|s_j-s_k|^{2\a_1}+ \min_{0 \le k \le j-1}|t_j-t_k|^{2\a_2}
\bigg).
\end{split}
\end{equation}

To prove the existence and continuity of the intersection local
times of $B^{\a_1}$ and $B^{\a_2}$, we will make use of the
following lemmas. Lemma \ref{Lem:Xiao08} is similar to Lemma 8.6 in
Xiao (2009) whose proof is elementary. Lemma \ref{Lem:Xiao97-ext}
and Lemma \ref{Lem:Xiao97} are extensions of Lemma 2.3 in Xiao
(1997) and will be useful for dealing with anisotropy of the
Gaussian random field $X_0$. Lemma \ref{Lem:CD82}, due to Cuzick and
DuPreez (1982), is a technical lemma.

\begin{lemma}\label{Lem:Xiao08}
Let $\beta$, $\gamma$ and $p$ be
positive constants, then for all $A \in (0, 1)$
\begin{equation}
\int_0^1\frac{r^{p-1}}{\big(A+r^{\ga}\big)^\beta}\, dr
\asymp\left\{\begin{array}{ll}
A^{\frac p \ga -\beta} \qquad\qquad \qquad &\hbox{ if } \,\,\beta\gamma> p, \\
\log\big(1+A^{-1/\ga}\big)  &\hbox{ if }\,\,\beta\gamma=p, \\
1   &\hbox{ if } \,\, \beta\gamma<p.
\end{array}
\right.
\end{equation}
In the above,  $f(A) \asymp g (A)$  means that the ratio $f(A)/g(A)$
is bounded from below and above by positive constants that do not
depend on $A \in (0, 1)$.
\end{lemma}

\begin{proof}\
This can be verified directly and we omit the details.
\end{proof}

\begin{lemma}\label{Lem:Xiao97-ext}
Let $\beta$, $\gamma$ and $p$ be positive constants such that
$\gamma\beta  \ge p$.
\begin{itemize}
\item[(i).] If $\gamma\beta  > p$,  then there exists
a constant $c_{_{2,3}}>0$ whose value depends on $\ga,\,\beta$ and
$p$ only such that for all $A \in (0, 1)$, $r >0$, $u^*\in\R^p$,
all integers $n\ge 1$ and all distinct $u_1,\ldots,u_n\in
O_p(u^*,r)$ we have
\begin{equation}\label{Eq:factorial2}
\int_{O_p(u^*,r)}\frac{du}{\left[ A+ \min\{|u-u_j|^{\ga},\,
j=1,\ldots,n\}\right]^{\beta}} \le c_{_{2,3}}\, n\,  A^{\frac{p}
{\ga}-\beta}.
\end{equation}

\item[(ii).] If $\gamma\beta = p$,  then for any $\kappa \in (0, 1)$ there exists
a constant $c_{_{2,4}}>0$ whose value depends on $\ga,\,\beta, \,
\kappa$ and $p$ only such that for all $A \in (0, 1)$, $r >0$,
$u^*\in\R^p$, all integers $n\ge 1$ and all distinct
$u_1,\ldots,u_n\in O_p(u^*,r)$ we have
\begin{equation}\label{Eq:factorial2-b}
\int_{O_p(u^*,r)}\frac{du}{\left[ A+ \min\{|u-u_j|^{\ga},\,
j=1,\ldots,n\}\right]^{\beta}} \le c_{_{2,4}}\, n\,\log \bigg[ e +
\Big(A^{-1/\ga} \frac r {n^{1/p}}\Big)^\kappa\bigg].
\end{equation}
\end{itemize}
\end{lemma}

\begin{proof}\ The idea of proof is similar to that of Lemma 2.3 in
Xiao (1997). Let
$$
\Gamma_i = \biggl\{ u \in  O_p(u^*,r) :\ |u - u_i| = \min\{|u -
u_j|, \ j = 1, \cdots , n\}\biggr\}.$$ Then
\begin{equation}\label{Eq:O1}
O_p(u^*,r) = \bigcup_{i=1}^n \Gamma_i \ \ \hbox{and}\ \ \
\la_p(O_p(u^*,r))= \sum_{i=1}^n \la_p(\Gamma_i).
\end{equation}
For every $u \in \Gamma_i$, we write $u$ $ = u_i + \rho \theta$, where
$0 \le \rho \le \rho_i(\theta)$ and $\theta \in S_{p-1}$, the unit
sphere in $\R^p$. Then
\begin{equation}\label{Eq:O2}
\begin{split}
\la_p (\Ga_i) &= C_p \int_{S_{p-1}} \nu(d
\theta)\int_0^{\rho_i(\theta)} \rho^{p-1}
d \rho\\
&= \frac{C_p} p\int_{S_{p-1}} \rho_i(\theta)^{p} \nu(d \theta),
\end{split}
\end{equation}
where $\nu$ is the normalized surface area in $S_{p-1}$ and $C_p$ is
a positive constant depending on $p$ only.

Denote the integral in (\ref{Eq:factorial2}) and (\ref{Eq:factorial2-b})
by $I_1$. We first consider the case of $\gamma\beta > p$.  By (\ref{Eq:O1}),
a change of variables and Lemma \ref{Lem:Xiao08}, we can write $I_1$ as
\begin{equation}
\begin{split}
I_1 &= \sum_{i=1}^n  \int_{\Ga_i} \frac{du} {\left [A +
\min\{|u-u_j|^{\ga},\,
j=1,\ldots,n\}\right]^{\beta} }\\
&= \sum_{i=1}^n C_p \int_{S_{p-1}} \nu(d \theta)
\int_0^{\rho_i(\theta)} \frac{
\rho^{p-1}} {\left(A +  \rho^{\ga}\right)^\beta}\,  d \rho\\
&= \sum_{i=1}^n C_p \,A^{\frac{p} {\ga}-\beta}\, \int_{S_{p-1}}
\nu(d \theta)
  \int_0^{A^{-1/\ga}\rho_i(\theta)} \frac{ \rho^{p-1}}
  {(1+  \rho^\ga)^\beta  }\, d \rho  \\
  &\le c_{_{2,3}}\, \sum_{i=1}^n A^{\frac{p} {\ga}-\beta} \int_{S_{p-1}}  \nu(d \theta)\\
  &= c_{_{2,3}}\, n\, A^{\frac{p} {\ga}-\beta}.
\end{split}
\end{equation}
This proves inequality (\ref{Eq:factorial2}).

Now we assume $\gamma\beta = p$. As above, we use (\ref{Eq:O1}) and
a change of variables to get
\begin{equation}
\begin{split}
I_1 &=\sum_{i=1}^n C_p \int_{S_{p-1}} \nu(d \theta)
\int_0^{A^{-1/\ga}\rho_i(\theta)}
\frac{\rho^{p-1}} {(1+\rho^\ga)^\beta}\,d \rho\\
&\le \frac{2 C_p} {\kappa}\, \sum_{i=1}^n \int_{S_{p-1}}\log\bigg[ e
+ \Big(A^{-1/\ga} \rho_i(\theta)\Big)^\kappa \bigg]\,\nu(d \theta).
\end{split}
\end{equation}
In the above, we have used the fact that if $\gamma\beta = p$ and
$\kappa \in (0, 1)$, then for all $x \ge 0$
$$\int_0^x \frac{\rho^{p-1}} {(1+\rho^\ga)^\beta}\,d \rho
\le \frac 2 \kappa \log(e + x^{\kappa}).$$

Since the function $\psi_1(x) = \log (e + x^{\kappa/p})$ is concave
on $(0, \infty)$, we apply (\ref{Eq:O2}) and Jensen's inequality
twice to derive
\begin{equation}
\begin{split}
I_1 &\le c_{_{2,5}} \,n \sum_{i=1}^n \frac 1 n\, \psi_1\Big( A^{-p/\ga}\la_p(\Ga_i) \Big)\\
&\le c_{_{2,4}}\, n \, \log \bigg[e + \Big(A^{-1/\ga} \frac r
{n^{1/p}}\Big)^\kappa\bigg].
\end{split}
\end{equation}
This finishes the proof of (\ref{Eq:factorial2-b}).
\end{proof}

\begin{lemma}\label{Lem:Xiao97}
Let $\beta>0$ be a constant and let $p\ge 1$ be an integer such that
$\beta< p$. Then the following statements hold:
\begin{itemize}
\item [(i).] For all $r >0$, $u^*\in\R^p$,
all integers $n\ge 1$, and all distinct $u_1,\ldots,u_n \in
O_p(u^*,r)$, we have
\begin{equation}\label{Eq:factorial3}
\int_{O_p(u^*,r)}\frac{du}{\min\{|u-u_j|^{\beta},\,j=1, \ldots,n\}}
\le c_{_{2,6}}\, n^{\frac{\beta} p} \, r^{p-\beta},
\end{equation}
where $c_{_{2,6}}>0$ is a constant whose value depends on $\beta$
and $p$ only.
\item[(ii).] For all constants $r >0$ and $K > 0$, all
$u^*\in\R^p$, integers $n\ge 1$, and all distinct $u_1,\ldots,u_n \in
O_p(u^*,r)$, we have
\begin{equation}\label{Eq:factorial4}
\begin{split}
&\int_{O_p(u^*,r)}  \log \bigg[e + K \big(\min\{|u-u_j|,\,j=1,
\ldots,n\}\big)^{-\beta} \bigg]\, du\\
&\qquad  \le c_{_{2,7}}\,r^p\, \log \bigg[e+K\Big(\frac r
{n^{1/p}}\Big)^{-\beta}\bigg],
\end{split}
\end{equation}
where $c_{_{2,7}}>0$ is a constant whose value depends on $\beta$
and $p$ only.
\end{itemize}
\end{lemma}

\begin{proof}\ Part (i) is a special case of Lemma 2.3 in
Xiao (1997). Hence, it only remains to prove Part (ii).
Denote the integral in (\ref{Eq:factorial4}) by $I_2$. As
in the proof of Lemma \ref{Lem:Xiao97-ext}, we have
\begin{equation}\label{Eq:factorial5}
\begin{split}
I_2= &\sum_{i=1}^n  \int_{\Ga_i} \log \Big[e+ K\,
\big(\min\{|u-u_j|,\,j=1,
\ldots,n\} \big)^{-\beta}\Big]\, du \\
&= \sum_{i=1}^n C_p \int_{S_{p-1}} \nu(d \theta)
\int_0^{\rho_i(\theta)}
\rho^{p-1}  \log\big( e+ K\,  \rho^{-\beta}\big)\, d \rho\\
&= \sum_{i=1}^n C_p \int_{S_{p-1}}\rho_i(\theta)^{p}\, \nu(d \theta)
\int_0^1 \rho^{p-1} \log\big(e+ K\rho_i(\theta)^{-\beta}
\rho^{-\beta}
\big)\, d \rho  \\
&\le  c_{_{2,8}}\,\sum_{i=1}^n \int_{S_{p-1}}\rho_i(\theta)^{p}
\log\big(e  +K\rho_i(\theta)^{-\beta} \big)\, \nu(d \theta).
\end{split}
\end{equation}
In deriving the last inequality, we have use the fact that $\log (e
+ x y) \le \log (e +x) + \log (e + y)$ for all $x, y \ge 0$. Since
$\beta< p$, we can verify that the function $ \psi_2(x) = x \log (e
+ K\,x^{-\beta/p}) $ is concave on $(0, \infty)$. By using Jensen's
inequality twice, we obtain
\begin{equation}\label{Eq:factorial6}
\begin{split}
I_2 &\le c_{_{2,8}}\, \sum_{i=1}^n  \psi_2\bigg(\int_{S_{p-1}}
\rho_i(\theta)^{p}\, \nu(d \theta)\bigg)\\
&\le c_{_{2,8}}\, n \,\psi_2\bigg(\frac 1 n \sum_{i=1}^n \la_p(\Gamma_\ell)\bigg)
\le c_{_{2,7}}\, n\, \psi_2\Big(\frac{r^p} n\Big).
\end{split}
\end{equation}
This finishes the proof of Lemma \ref{Lem:Xiao97}.
\end{proof}

\begin{lemma} \label{Lem:CD82}
Let $Z_1, \ldots, Z_n$ be the mean zero Gaussian random variables
which are linearly independent and assume that
$$\int_{- \infty}^{\infty} g(v) e^{- \ep v^2} dv < \infty$$
for all $\ep > 0$. Then
\begin{equation}
\begin{split}
&\int_{{\R}^n} g(v_1) \exp \Bigl[- \frac1 2 {\rm
Var}\Big(\sum_{j=1}^n v_j Z_j\Big)
\Bigr] dv_1 \cdots dv_n \\
&= \frac{(2 \pi)^{(n-1)/2}} {({\rm detCov}(Z_1, \ldots,
Z_n))^{1/2}}\ \int_{- \infty}^{\infty} g\Bigl(\frac {v} {
\si_1}\Bigr) e^{- v^2/2}\, dv,
\end{split}
\end{equation}
where $\si_1^2 = {\rm Var}(Z_1|Z_2, \ldots, Z_n)$ is the conditional
variance of $Z_1$ given $Z_2,$ $ \ldots,  Z_n$.
\end{lemma}

\section{Intersection local times and their joint continuity}
\label{Sec:LT}

In this section, we briefly recall the definition of local time
as occupation density [cf. Geman and Horowitz (1980)] and then study the existence and joint continuity
of the intersection local times of $B^{\a_1}$ and $B^{\a_2}$.

Let $Y(t)$ be a [random] Borel vector field on $\R^p$ with values in $\R^q$.
For any Borel set $E \subseteq \R^p$, the occupation measure of $Y$
on $E$ is defined as the following measure on $\R^q$:
\[
   \mu_{_E}(\bullet) = \la_p \big\{ t \in E: Y(t) \in \bullet
   \big\}.
\]

If $\mu_{_E}$ is absolutely continuous with respect to the Lebesgue
measure $\la_q$, we say that $Y(t)$ has \emph{local time} on $E$,
and define its local time, $L(\bullet, E)$, as the Radon--Nikod\'ym
derivative of $\mu_{_E}$ with respect to $\lambda_q$, i.e.,
\[
    L(x,E) = \frac{d\mu_{_E}} {d\lambda_q}(x),\qquad \forall x\in\R^q.
\]
In the above, $x$ is the so-called \emph{space variable}, and $E$ is
the \emph{time variable}. Note that if $Y$ has local times on $E$
then for every Borel set $F \subseteq E$, $L(x, F)$ also exists.

It follows from Theorem 6.4 in Geman and Horowitz (1980) that the
local time has a measurable modification that satisfies the
following \emph{occupation density formula}: For every Borel set $E
\subseteq \R^p$, and for every measurable function $f : \R^q \to
\R_+$,
\begin{equation}\label{Eq:occupation}
\int_E f(Y(t))\, d t = \int_{\R^q} f(x) L(x, E)\, dx.
\end{equation}

Suppose we fix a rectangle $E = [a,\, a+h] \subseteq \R^p$, where $a
\in \R^p$ and $h \in \R_+^p$. If we can choose a version of the
local time, still denoted by $L(x, [a,\, a+t])$,  such that it is a
continuous function of $(x, t)$ $\in$ $\R^q\times [0,h]$,  $Y$ is
said to have a \emph{jointly continuous local time} on $E$. When a
local time is jointly continuous, $L(x, \cdot)$ can be extended to
be a finite Borel measure supported on the level set
\begin{equation}\label{eq:inv}
    Y_E^{-1}(x) = \{ t \in E: Y(t) = x\};
\end{equation}
see Theorem 8.6.1 in Adler (1981) for details. This makes local
times, besides of interest on their own right, a useful tool in
studying fractal  properties of $Y$.

It follows from (25.5) and (25.7) in Geman and Horowitz (1980) that,
for all $x, y \in \R^q$, $E\subseteq \R^p$ a closed interval and all
integers $n \ge 1$,
\begin{equation} \label{Eq:moment1}
\begin{split}
\E \Big[ L(x, E)^n\Big] &= (2\pi)^{-nq} \int_{E^n} \int_{{
\R}^{nq}} \exp \bigg(- i \sum_{j=1}^{n} \l u_j, x\r \bigg)\\
&\qquad \qquad \qquad \qquad \ \ \times \E \exp\bigg( i \sum_{j=1}^{n} \l
u_j, Y(t_j)\r \bigg)\, d\overline{u}\ d\overline{t}
\end{split}
\end{equation}
and, for all even integers $n \ge 2$,
\begin{equation} \label{Eq:moment2}
\begin{split}
\E\Big[ (L(x, E) - L(y, E))^n\Big] =  &(2 \pi)^{- nq} \int_{E^n}
\int_{{\R}^{nq}} \prod_{j=1}^{n} \Big[ e^{- i \l u_j, x\r} - e^{-
i \l u_j, y \r}\Big]\\
& \qquad \qquad \times \E \exp \bigg( i \sum_{j=1}^n \l u_j, Y(t_j)
\r\bigg) d\overline{u}\ d\overline{t},
\end{split}
\end{equation}
where $\overline{u} = ( u_1, \ldots, u_n),\ \overline{t} = (t_1,
\ldots, t_n),$ and each $u_j \in \R^q,\ t_j \in E.$ In the
coordinate notation we then write  $u_j = (u_{j,1}, \ldots,
u_{j,q}).$

The main results of this section are the following Theorem \ref{Thm:Exist} and
Theorem \ref{Thm:JCont} for the existence and the joint continuity of the
intersection local times of two independent fractional Brownian
motions in $\R^d$.

\begin{theorem}\label{Thm:Exist}
Let $B^{\a_1}= \{B^{\a_1}(s),\,s\in\R^{N_1}\}$ and
$B^{\a_2}= \{B^{\a_2}(t),\,t\in\R^{N_2}\}$ be two independent fractional
Brownian motions with values in $\R^d$ and Hurst indices $\a_1$ and $\a_2$,
respectively. Let $X=\{X(s,t), (s,t) \in \R^N\}$ be the $(N,d)$-Gaussian random
field defined by (\ref{Def:X}). Then, for any given constant $R>0$,
$X$ has a local time $L\big(x, O_{N_1,N_2}(0,R)\big)\in
L^2(\P\times\la_d)$ if and only if
$\frac{N_1}{\a_1}+\frac{N_2}{\a_2}>d$. Furthermore, if it exists,
the local time of $X$ admits the following $L^2$-representation
\begin{equation}\label{Eq:Int_LT}
L\big(x, O_{N_1,N_2}(0,R)\big)=(2\pi)^{-d}\int_{\R^d} e^{-i \l y,
x\r} \int_{O^2_{N_1,N_2}(0,R)}e^{i\l y, B^{\a_1}(s)-B^{\a_2}(t)\r}\,
dsdtdy,
\end{equation}
and the local time $L$ can be chosen as a kernel $L\big(\cdot,\, \cdot\big)$
on $\R^d \times {\cal B}\big( O_{N_1,N_2}(0,R)\big)$.
In particular, if $\frac{N_1}{\a_1}+\frac{N_2}{\a_2}>d$, then
$B^{\a_1}$ and $B^{\a_2}$ have an intersection local time which can
be defined as $ L^{\a_1,\a_2}\big(O_{N_1,N_2}(0,R)\big) := L \big(0,
O_{N_1,N_2}(0,R)\big).$
\end{theorem}

Some remarks about Theorem \ref{Thm:Exist} are in order.

\begin{remark}\
(i) When $N_1=N_2=1$, $\a_1 = \a_2 = H$ and $H d <2$,
the existence of the intersection local time was proved by Nualart
and Oritz-Latorre (2007) as the $L^2$-limit of
\begin{equation}\label{Eq:I_eps}
I_\ep(B^{H},\widetilde{B}^H) \equiv \int_0^R\int_0^R
p_\ep\big(B^{H}(s) - \widetilde{B}^{H}(t)\big)\, ds dt, \quad \hbox{
as } \ \ep \to 0,
\end{equation}
where $p_\ep(x) = (2\pi \ep)^{d/2} \exp(- |x|^2/(2\ep))$. They also
proved that if $H d \ge 2$, then
\[
\lim_{\ep \downarrow 0} \E\big[I_\ep(B^{H},\widetilde{B}^H)\big] =
\infty \ \ \hbox{ and }\ \ \lim_{\ep \downarrow 0} {\rm
Var}\big[I_\ep(B^{H},\widetilde{B}^H)\big] = \infty.
\]
In the above, $B^{H} = \{B^H(t), t \ge 0\}$ and $\widetilde{B}^H =
\{\widetilde{B}^H(t), t \ge 0\}$ are two independent fractional
Brownian motions with values in $\R^d$ and index $H \in (0, 1)$.
Similar method can be applied to show that the intersection local time
$L^{\a_1,\a_2}\big(O_{N_1,N_2}(0,R)\big)$ in Theorem \ref{Thm:Exist} can be
chosen as the $L^2$-limit of the following approximating functionals
\begin{equation}\label{Eq:I_eps-2}
I_\ep\big(B^{\a_1},{B}^{\a_2}\big) \equiv \int_{O_{N_1,N_2}(0,R)}
p_\ep\big(B^{\a_1}(s) - B^{\a_2}(t)\big)\, ds dt, \quad \hbox{
as } \ \ep \to 0.
\end{equation}
Moreover,
we are able to show that, if $\frac{N_1}{\a_1}+\frac{N_2}{\a_2} =d$,
then
\begin{equation}\label{Eq:e0}
\E\big[I_\eps(B^{\a_1},B^{\a_2})\big]\sim
c(\a_1,\a_2,N_1,N_2)\ln\left(\frac1\eps\right), \quad \hbox{ as } \
\ep \to 0,
\end{equation}
where $c(\a_1,\a_2,N_1,N_2)>0$ is a constant depending on
$\a_1,\,\a_2$ and $N_1,\,N_2$ only. This raises an interesting
question whether $I_\ep$ can be renormalized to converge to a
non-trivial limiting process. This and other related questions will
be dealt with elsewhere since they require different methods.

(ii) It follows from the operator-scaling property (\ref{Eq:OSS}) of $X$ and
(\ref{Eq:Int_LT}) that the intersection local time $ L^{\a_1,\a_2}
\big(O_{N_1,N_2}(0,R)\big)$ has the following scaling property: For
any constant $c > 0$,
\begin{equation}\label{Eq:Lss}
\Big\{L^{\a_1,\a_2} \big(c^A\,O_{N_1,N_2}(0,R)\big), \, R>0\Big\}
\stackrel{d} {=} \Big\{c^{\frac{N_1}{\a_1} + \frac{N_2}{\a_2} - d}\,
L^{\a_1,\a_2} \big(O_{N_1,N_2}(0,R)\big), \, R>0\Big\}.
\end{equation}
Here $A$ is the $N\times N$ diagonal matrix as in (\ref{Eq:OSS}).

(iii) We say that the sample functions of $B^{\a_1}$ and $B^{\a_2}$
intersect if there exist $s \in \R^{N_1}$ and $t\in \R^{N_2}$ such
that $B^{\a_1} (s)= B^{\a_2}(t)$. It is also of interest to study
the geometric properties of the set of intersection times
\[
M_2= \{(s, t) \in \R^{N}:  B^{\a_1} (s)= B^{\a_2}(t)\}
\]
and the set of intersection points
\[
D_2 = \{x \in \R^d:  x =B^{\a_1} (s)= B^{\a_2}(t) \ \hbox{ for some
}\, (s, t) \in \R^N\},
\]
because they are often random fractals. The existence of the
intersection local time and its properties are closely related to
the existence of intersections of the sample functions of $B^{\a_1}$
and $B^{\a_2}$ and the geometric properties of $M_2$ and $D_2$.
Similar to Theorem 7.1 in Xiao (2009), we can prove that if
$\frac{N_1}{\a_1}+\frac{N_2}{\a_2} >d$ then $M_2 \ne \emptyset$ with
positive probability. On the other hand, Theorem 3.2 in
Xiao (1999) proved that if  $\frac{N_1}{\a_1}+\frac{N_2}{\a_2}\le
d$ then $M_2 = \emptyset$ almost surely. In Section \ref{Sec:H}, we will give
more information on the Hausdorff and packing dimensions of $M_2$, as well
as a lower bound for the exact Hausdorff measure of $M_2$.  The
Hausdorff and packing dimensions of $D_2$ are determined in Section
\ref{Sec:dim}.
\end{remark}

\begin{proof}{\bf of Theorem \ref{Thm:Exist}} Note that the
Fourier transform of the
occupation measure $\mu_{_{O_{N_1,N_2}(0,R)}}$ of $X$ is
\[
\widehat{\mu}_{_{O_{N_1,N_2}(0,R)}}(\xi)=\int_{O_{N_1,N_2}(0,R)}
e^{i\l\xi,X(s,t)\r}dsdt.
\]
It follows from the Plancherel Theorem that $X$ has a local time
$L\big(x,O_{N_1,N_2}(0,R)\big)\in L^2(\P\times\la_d)$ with a
representation (\ref{Eq:Int_LT}) if and only if
\begin{equation}\label{Eq:Exist1}
{\cal
J}:=\int_{O^2_{N_1,N_2}(0,R)}dsdtdvdw\int_{\R^d}\left|\E\exp\Big(i\l
y, X(s,t)-X(v,w)\r\Big)\right|dy<\infty.
\end{equation}
See Theorem 21.9 of Geman and Horowitz (1980). Hence, it suffices
to prove that Eq. (\ref{Eq:Exist1}) holds if and only if
$\frac{N_1}{\a_1}+\frac{N_2}{\a_2}>d$. For this purpose, we use the
independence of the coordinate processes of $X$, (\ref{Def:X})
and (\ref{Eq:Cov_fBm}) to deduce that
\begin{equation}\label{Eq:Exist2}
\begin{split}
{\cal J} &=\int_{O^2_{N_1,N_2}(0,R)}\frac{dsdtdvdw}
{\left[\E\big(X_0(s,t)-X_0(v,w)\big)^2\right]^{d/2}}\\
&=\int_{O^2_{N_1,N_2}(0,R)}\frac{dsdtdvdw}{\left[|s-v|^{2\a_1}
+|t-w|^{2\a_2}\right]^{d/2}}.
\end{split}
\end{equation}
By using spherical variable substitutions and Lemma \ref{Lem:Xiao08},
it is elementary to verify that the last integral in Eq. (\ref{Eq:Exist2})
is finite if and only if $\frac{N_1}{\a_1}+\frac{N_2}{\a_2}>d$.

When the later holds, one can apply Theorem 6.3 in Geman and Horowitz (1980)
to choose a version of the local time of $X$, still denoted by $L$, such
that it is a kernel in the following sense:
For every $x \in \R^d$, $L(x, \cdot)$ is a finite measure on ${\cal
B}(O_{N_1,N_2}(0,R))$ and, for every Borel set $E\in {\cal
B}(O_{N_1,N_2}(0,R))$, $x \mapsto L(x, E)$ is a
measurable function. This proves the main conclusion of Theorem
\ref{Thm:Exist}. Finally, by taking $x = 0$ we prove the last conclusion of
Theorem \ref{Thm:Exist}.
\end{proof}

\begin{theorem} \label{Thm:JCont}
Let $B^{\a_1}$ and $B^{\a_2}$ be defined as that in Theorem \ref{Thm:Exist}. If
$\frac{N_1}{\a_1}+\frac{N_2}{\a_2}>d$, then $B^{\a_1}$ and
$B^{\a_2}$ have almost surely a continuous intersection local time
on $\R^{N_1 + N_2}$.
\end{theorem}

As in the proof of Theorem \ref{Thm:Exist}, we will prove a stronger
result that $X$ has almost surely a jointly continuous local time on
$\R^{N_1 + N_2}$. The proof is based on the following Lemma
\ref{Lem:MomEst} and Lemma \ref{Lem:MomEst2}. They will also play
an essential r\^{o}le in Section \ref{Sec:H} for establishing
H\"{o}lder conditions for the intersection local times.

Under the condition $\frac{N_1}{\a_1}+\frac{N_2}{\a_2}>d$, define
\begin{equation}\label{Def:tau}
\tau=\left\{\begin{array}{ll}
1 \qquad &\hbox{ if }\,\frac{N_1}{\a_1}>d, \\
2 \qquad &\hbox{ if }\,\frac{N_1}{\a_1}\le
d<\frac{N_1}{\a_1}+\frac{N_2}{\a_2}
            \end{array}
\right.
\end{equation}
and
\begin{equation}\label{Def:beta}
\beta_\tau = \left\{\begin{array}{ll}
N-\a_1d \qquad \quad &\hbox{ if } \, \tau = 1,\\
N_2+\frac{\a_2}{\a_1}N_1-\a_2d &\hbox{ if }\, \tau = 2.\\
\end{array}
\right.
\end{equation}
[Recall that we assumed $\a_1\le \a_2$ throughout the paper, and
$N=N_1+N_2$.] We will also make use of the following notation:
\begin{equation}\label{Def:gamma}
\eta_\tau = \left\{\begin{array}{ll}
\frac{\a_1d} {N_1} \qquad \quad &\hbox{ if } \, \tau = 1,\\
\frac{\a_2 d} {N_2} + 1 - \frac{\a_2 N_1} {\a_1 N_2}  &\hbox{ if }\, \tau = 2.\\
\end{array}
\right.
\end{equation}
Note that, if $N_1 = \a_1 d$, then $\beta_\tau = N_2$ and $\eta_\tau
= 1$. To emphasize the importance of $\beta_\tau$ and $\eta_\tau$,
we point out that $\beta_\tau$ is the Hausdorff dimension of the set $M_2$
of intersection times and $\eta_\tau$ is useful for determining the exact Hausdorff
measure of $M_2$. See Section \ref{Sec:H} for more information.

\begin{lemma}\label{Lem:MomEst}
Suppose the assumptions of
Theorem \ref{Thm:JCont} hold. Then, there exist positive and finite constants
$\eps \in (0, 1/e)$ and $c_{_{3,1}}$, which depend on
$\a_1,\a_2,\,N_1,\,N_2$ and $d$ only, such that for all
$r\in(0,\,\eps),$  $D:=O_{N_1,N_2}(u,r)$, where $u = (u_1, u_2) \in
\R^N$,
 all $x\in\R^d$
and all integers $n\ge 1$, we have
\begin{equation}\label{Eq:momentEst1}
\E\left[L(x,D)^n\right]\le \left\{\begin{array}{ll} c_{_{3,1}}^n\,
(n!)^{\eta_1}\,
r^{n\,\beta_1} \quad &\hbox{ if } \frac{N_1}{\a_1}>d,\\
c_{_{3,1}}^n\, n!\, r^{n\,N_2} \, \prod_{j=1}^n\log \Big(e+\frac
{j^{\left(\frac{\a_2} {N_2} - \frac {\a_1} {N_1}\right)^+ }}
{r^{\a_2- \a_1 }} \Big)\quad &\hbox{ if } \frac{N_1}{\a_1}=d,\\
c_{_{3,1}}^n\, (n!)^{\eta_2}\, r^{n\,\beta_2} \quad &\hbox{ if }
\frac{N_1}{\a_1}< d < \frac{N_1}{\a_1}+\frac{N_2}{\a_2}.
\end{array}\right.
\end{equation}
In the above, $y^+ = \max \{y, 0\}$ for every $y\in\R$.
\end{lemma}

\begin{remark}\ From (\ref{Def:beta}) and (\ref{Def:gamma}), it can be
verified that
\begin{equation}\label{Eq:betaga}
\frac {N - \beta_\tau} N \le \eta_\tau \le N- \beta_\tau.
\end{equation}
We observe that the power of $n!$ in (\ref{Eq:momentEst1}) becomes $(N -
\beta_\tau)/N$ when $X$ is an isotropic Gaussian field as in Xiao (1997)
and is $N- \beta_\tau$ when $X$ is anisotropic in every coordinate (with the
same scaling or H\"older index) as in Ayache, Wu and Xiao (2008).
These seem to be the extreme cases. In the present paper, if we
assume  $N_1 \ne N_2$ and $\a_1 \ne \a_2$, then strict inequalities
in \eqref{Eq:betaga} may hold and if, in addition, $N_1 = \a_1 d$,
then extra logarithmic factors appear in the estimate
(\ref{Eq:momentEst1}). Lemma \ref{Lem:MomEst} suggests that the local time $L(x,
\cdot)$ may satisfy a law of the iterated logarithm which is
different from those for the local times of an $(N, d)$-fractional
Brownian motion or an $(N, d)$-fractional Brownian sheet with index
$(\a, \ldots, \a)$; see (\ref{Eq:Holder0}), (\ref{Eq:Holder-fbm})
and (\ref{Eq:Holder-sh}) below. This leads us to expect that the exact
Hausdorff measure function for $M_2$ may be different from those for
the level sets of fractional Brownian motion and fractional Brownian
sheets, respectively. It would be interesting to investigate these
problems.
\end{remark}

\begin{proof}{\bf of Lemma \ref{Lem:MomEst}} Even though the
proof of Lemma \ref{Lem:MomEst} follows
the same spirit of the proofs of Lemma 2.5 in Xiao (1997) and
Lemma 3.7 in Ayache, Wu and Xiao (2008), there are
some subtle differences [see the remark above]. Hence
we give a complete proof. In particular, we
provide a direct way to estimate the last integral in (\ref{Eq:m51})
below. We believe that this method will be useful elsewhere.

It follows from (\ref{Eq:moment1}) and the fact that $X_1, \ldots,
X_d$ are independent copies of $X_0$ that, for all integers $n \ge
1$,
\begin{equation}\label{Eq:m51}
\begin{split}
\E \big[ L(x, D)^n\big] &\le (2 \pi)^{- nd} \int_{D^n}\prod_{k=1}^d
\bigg\{\int_{{\R}^n}  \exp\bigg[ -\frac1 2 {\rm Var} \bigg(
\sum_{j=1}^n u_{j,k}\, X_0(s_j,t_j)\bigg) \bigg]\, d \overline{u}_k
\bigg\}\,
d\overline{\t}\\
& = (2 \pi)^{- nd/2} \int_{D^n}  \Big[{\rm
detCov}\left(X_0(s_1,t_1), \ldots, X_0(s_n,t_n)\right)\Big]^{-\frac
d 2}\, d\overline{\t},
\end{split}
\end{equation}
where $\overline{u}_k = (u_{1,k}, \ldots, u_{n,k}) \in {\R}^n$,
$\overline{\t} =(s_1,t_1,\ldots,s_n,t_n)$ and the equality follows
from the fact that for any positive definite $n \times n$ matrix
$\Gamma$,
\begin{equation}\label{Eq:Normal}
\int_{\R^n} \frac {[{\rm det}(\Gamma)]^{1/2}} {(2\pi)^{n/2}}\,
\exp\Big(- \frac 1 2 x'\Gamma x\Big)\, dx = 1.
\end{equation}

In order to prove Eq. (\ref{Eq:momentEst1}), we consider the three
cases separately: \ $\frac{N_1}{\a_1}>d$, $\frac{N_1}{\a_1}<
d<\frac{N_1}{\a_1}+\frac{N_2}{\a_2}$ and $\frac{N_1}{\a_1}=d$.

In the case that $\frac{N_1}{\a_1}>d$, thanks to Eq.
(\ref{Eq:Det_XB}), we have
\begin{equation}\label{Eq:m520b}
\begin{split}
\E \big[ L(x, D)^n\big] &\le c_{_{3,2}}^n\int_{D^n}
\prod_{j=1}^n\frac1{\min\{|s_j-s_i|^{\a_1 d},\,\,0\le
i\le j-1\}}\, d\overline{s}\, d\overline{t}\\
& = c_{_{3,2}}^n\int_{O^n_{N_2}(u_2,r)}\left(\int_{O^n_{N_1}(u_1,r)}
\prod_{j=1}^n\frac1{\left(\min\{|s_j-s_i|^{\a_1},\,\, 0\le i\le
j-1\}\right)^{d }}\, d\overline{s}\right) d\overline{t}\\
&  = c_{_{3,3}}^nr^{n N_2} \, \int_{O^n_{N_1}(u_1,r)}
\prod_{j=1}^n\frac1{\min\{|s_j-s_i|^{\a_1 d},\,\, 0\le i\le j-1\}}\,
d\overline{s},
\end{split}
\end{equation}
where $s_0:=0$, $\overline{s}=(s_1,\ldots,s_n)$ and
$\overline{t}=(t_1,\ldots,t_n)$.

Since $N_1>\a_1 d$, we integrate the last integral in Eq. (\ref{Eq:m520b})
in the order $ds_n, \ldots, ds_1$ and apply Part (i) of Lemma \ref{Lem:Xiao97}
iteratively. This yields
\begin{equation}\label{Eq:m520c}
\begin{split}
\E \big[ L(x, D)^n\big]& \le c_{_{3,1}}^n (n!)^{\frac{\a_1
d}{N_1}}\, r^{n (N_1-\a_1 d)}\times r^{nN_2}  = c_{_{3,1}}^n
(n!)^{\eta_1}\,r^{n\beta_1},
\end{split}
\end{equation}
which proves Eq. (\ref{Eq:momentEst1}) for the case
$\frac{N_1}{\a_1}>d$ [i.e., $\tau=1$].

In the second and third cases [i.e., $\frac{N_1}{\a_1}\le
d<\frac{N_1}{\a_1}+\frac{N_2}{\a_2}$] we use (\ref{Eq:m51}) and
(\ref{Eq:Det_XB}) to obtain
\begin{equation}\label{Eq:case2-1}
\E \big[ L(x, D)^n\big] \le c_{_{3,4}}^n \int_{D^n} \prod_{j=1}^n
\frac 1 {\big( \min\limits_{0 \le k \le j-1}|s_j-s_k|^{\a_1}+ \min\limits_{0 \le k
\le j-1}|t_j-t_k|^{\a_2}\big)^d}\, d\overline{s}d\overline{t}.
\end{equation}
To estimate the last integral in (\ref{Eq:case2-1}), we will
integrate in the order of $ds_n, dt_n,\ldots, ds_1, dt_1$. In the
case of $\frac{N_1}{\a_1} < d<\frac{N_1}{\a_1}+\frac{N_2}{\a_2}$, we
apply Part (i) of Lemma \ref{Lem:Xiao97-ext} with $A = \min_{0 \le k \le
n-1}|t_n-t_k|^{\a_2}$ to derive
\begin{equation}\label{Eq:case2-2}
\begin{split}
\int_{O_{N_1}(u_1,r)} \frac {ds_n} {\big( \min\limits_{0 \le k \le
n-1}|s_n-s_k|^{\a_1}+ \min\limits_{0 \le k \le n-1}|t_n-t_k|^{\a_2}\big)^d}
&\quad \le \frac{c_{_{2, 3}}\, n} {\big[\min\limits_{0 \le k \le
n-1}|t_n-t_k|^{\a_2}\big]^{d - \frac{N_1} {\a_1}}}.
\end{split}
\end{equation}
Since $\a_2 (d - \frac{N_1} {\a_1}) < N_2$, it follows from
(\ref{Eq:case2-2}) and Part (i) of Lemma \ref{Lem:Xiao97} that
\begin{equation}\label{Eq:case2-3}
\begin{split}
&\int_{D} \frac {ds_n dt_n} {\big( \min\limits_{0 \le k \le
n-1}|s_n-s_k|^{\a_1}+ \min\limits_{0 \le k \le n-1}|t_n-t_k|^{\a_2}\big)^d}\\
&\le c_{_{2, 3}}\, n\, \int_{O_{N_2}(u_2,r)} \frac{dt_n}
{\big(\min\limits_{0 \le k \le n-1}|t_n-t_k|^{\a_2}\big)^{d - \frac{N_1}
{\a_1}}} \\
&\quad \le c_{_{3,5}}\, n^{1 + \frac{\a_2(\a_1 d- N_1)} {\a_1
N_2}}\, r^{N_2 - \a_2 (d - \frac{N_1}{\a_1})} = c_{_{3,5}}\,
n^{\eta_2}\, r^{\beta_2}.
\end{split}
\end{equation}
Repeating the above procedure yields (\ref{Eq:momentEst1}) for the
case of $\frac{N_1}{\a_1} < d<\frac{N_1}{\a_1}+\frac{N_2}{\a_2}$.

Finally, we consider the case of $\frac{N_1}{\a_1} = d$. Let $\kappa
\in (0, 1)$ be a constant such that $\kappa \a_2/\a_1 < N_2$.
Applying Part (ii) of Lemma \ref{Lem:Xiao97-ext} with $A = \min_{0 \le k \le
n-1}|t_n-t_k|^{\a_2}$, we have
\begin{equation}\label{Eq:case2-4}
\begin{split}
&\int_{O_{N_1}(u_1,r)} \frac {ds_n} {\big(\min\limits_{0 \le k \le
n-1}|s_n-s_k|^{\a_1}+
\min\limits_{0 \le k \le n-1}|t_n-t_k|^{\a_2}\big)^d}\\
&\quad \le c_{_{3,6}}\, n \log \bigg[e + \bigg( \big(\min\limits_{0 \le k
\le n-1}|t_n-t_k|^{\a_2}\big)^{-1/\a_1} \, \frac r
{n^{1/N_1}}\bigg)^\kappa\bigg].
\end{split}
\end{equation}
It follows from (\ref{Eq:case2-4}) and Part (ii) of Lemma \ref{Lem:Xiao97} (with
$\beta = \kappa \a_2/\a_1$ and $K = (r\,
{n^{-1/N_1}})^\kappa$) that
\begin{equation}\label{Eq:case2-5}
\begin{split}
&\int_{D} \frac {ds_n dt_n} {\big( \min\limits_{0 \le k \le
n-1}|s_n-s_k|^{\a_1}+
\min\limits_{0 \le k \le n-1}|t_n-t_k|^{\a_2}\big)^d}\\
&\quad \le c_{_{3,6}}\, n\, \int_{O_{N_2}(u_2,r)}\log \bigg[e +
\big(\min_{0 \le k \le n-1}|t_n-t_k|\big)^{- \kappa \a_2/\a_1}
\, \Big(\frac r {n^{1/N_1}}\Big)^\kappa \bigg] \, dt_n\\
&\quad \le c_{_{3,7}}\, n \, r^{N_2}\, \log \bigg[e + \Big(\frac r
{n^{1/N_2}}\Big)^{-\kappa \a_2/\a_1} \Big(\frac r
{n^{1/N_1}}\Big)^\kappa
 \bigg] \\
 &\quad = c_{_{3,7}}n\, r^{N_2}\log\bigg[e+\Big(\frac{n^{\frac{\a_2}{N_2}-
 \frac{\a_1}{N_1}}}{r^{\a_2-\a_1}}\Big)^{\kappa/\a_1}\bigg]\\
 & \quad \le c_{_{3,8}}\, n\, r^{N_2} \,\log \bigg( e+\frac{n^{\left(\frac{\a_2}{N_2}-
 \frac{\a_1}{N_1}\right)^+}} {r^{\a_2-\a_1}}\bigg).
\end{split}
\end{equation}
[Recall that $y^+ = \max\{y, 0\}$.]

By iterating the procedure and integrating $ds_{n-1}, dt_{n-1},
\ldots, ds_1, dt_1$, we obtain that
\begin{equation}\label{Eq:m53}
\begin{split}
\E \big[ L(x, D)^n\big]\le c_{_{3,1}}^n\, (n!)^{ \eta_2}\, r^{n
\beta_2} \,\prod_{j=1}^n \log \bigg(e+\frac
{j^{\left(\frac{\a_2}{N_2}-
 \frac{\a_1}{N_1}\right)^+}} {r^{\a_2 -\a_1}} \bigg).
\end{split}
\end{equation}
This finishes the proof of the moment estimate
(\ref{Eq:momentEst1}).
\end{proof}
The following lemma estimates the higher moments of the increments of
the local times of $X$. Combined with Kolmogorov's continuity
theorem, it immediately implies the existence of a continuous
version of $x \mapsto L(x, D)$.

\begin{lemma}\label{Lem:MomEst2} Suppose the assumptions of
Theorem \ref{Thm:JCont} hold. Then, there exist positive constants $c_{_{3,9}}$
and $\kappa_1 $, depending on $\eps, \a_1,\a_2,\,N_1,\,N_2$ and $d$
only, such that, for any $r >0,$ $D:=O_{N_1,N_2}(u,r)$ for $u = (u_1,
u_2) \in \R^N$, all $x,\,y\in\R^d$ with $|x-y|\le 1$, all even
integers $n\ge1$, and all $\gamma\in(0,\,1)$ small enough, we have
\begin{equation}\label{Eq:momentEst2}
\E\left[\big(L(x,D)-L(y,D)\big)^n\right]\le
c_{_{3,9}}^n\,(n!)^{\eta_\tau +\kappa_1 \gamma}\, |x-y|^{n \gamma}\,
r^{n (\beta_\tau - \kappa_1 \gamma)}.
\end{equation}
\end{lemma}

\begin{proof}\  Let $\ga \in (0, 1)$ be a small constant whose value
will be determined later. Note that by the elementary inequalities
\begin{equation}
|e^{i u} - 1 | \le 2^{1 - \gamma} |u|^{\gamma} \qquad \hbox { for
all }\ u \in {\R}
\end{equation}
and $|u + v|^\gamma \le |u|^\gamma + |v|^\gamma$, we see that for
all $u_1, \ldots, u_n,\ x,\ y \in \R^d$,
\begin{equation}\label{Eq:gamma}
\prod_{j=1}^{n} \Big| e^{- i \l u_j, x\r} - e^{- i \l u_j, y \r
}\Big| \le 2^{(1 - \ga)n} \ |x-y|^{n \ga}\ {\sum}^{ '}
\prod_{j=1}^{n} |u_{j,k_j}|^ {\ga},
\end{equation}
where the summation $ \sum$\'\ is taken over all the sequences
$(k_1, \ldots, k_n) \in \{ 1, \ldots, d\}^n$.

It follows from (\ref{Eq:moment2}) and (\ref{Eq:gamma}) that  for
every even integer $n \ge 2$,
\begin{equation} \label{Eq:LTmoment3}
\begin{split}
\E\Big[ (L(x, &\,D) - L(y, D))^n\Big]  \le (2 \pi)^{-nd} 2^{(1 -
\gamma)n}\,     |x-y|^{n \gamma}  \\
&\qquad \times {\sum}^{'} \int_{D^n} \int_{{\R}^{n d} }
\prod_{m=1}^{n} |u_{m,k_m}|^{\gamma} \, \E \exp \bigg( - i
\sum_{j=1}^n \l u_j,\, X(s_j,t_j)\r  \bigg) \, d \overline {u}
\, d\overline{\t}  \\
& \qquad \le c_{_{3,10}}^n |x-y|^{n \gamma}\, {\sum}^{'}
\int_{D^n} \, d\overline{\t} \\
&\qquad  \times \prod_{m=1}^n \Bigg\{ \int_{\R^{nd}} |u_{m,k_m}|^{n
\gamma} \,  \exp\bigg[- \frac1 2 {\rm Var} \bigg(\sum^n_{j=1} \l
u_j,\, X(s_j,t_j)\r \bigg)\bigg]  \, d\overline{u}\Bigg\}^{1/n},
\end{split}
\end{equation}
where the last inequality follows from the generalized H\"older
inequality.

Now we fix a vector $\overline{k} = (k_1, k_2, \ldots, k_n) \in \{
1, \ldots, d\}^n$ and $n$ points $(s_1,t_1), \ldots, (s_n,t_n) \in
D\backslash\{0\} $ such that $s_1,\,t_1, \ldots, s_n,\,t_n$ are all
distinct [the set of such points has full $nN$-dimensional Lebesgue
measure]. Let $\mathcal{M} = \mathcal{M}( \overline{k},
\overline{\t}, \gamma)$ be defined  by
\begin{equation} \label{Eq:M1}
\mathcal{M} = \prod_{m=1}^n \Bigg\{ \int_{\R^{nd}} |u_{m,k_m}|^{n
\gamma} \, \exp\bigg[- \frac1 2 {\rm Var} \bigg(\sum^n_{j=1} \l
u_j,\, X(s_j,t_j)\r    \bigg)\bigg]  \, d\overline{u}\Bigg\}^{1/n}.
\end{equation}
Note that $X_\ell$ ($ 1 \le \ell \le d$) are independent copies of
$X_0$. By the strong local nondeterminism of fractional Brownian
motions $B_0^{\a_1}$ and $B_0^{\a_2}$ and Eq. (\ref{Eq:Det_XB}), the
random variables $X_\ell(s_j,t_j)$ ($ 1 \le \ell \le d, 1 \le j \le
n$) are linearly independent. Hence Lemma \ref{Lem:CD82} gives
\begin{equation}\label{Eq:M4}
\begin{split}
&\int_{\R^{nd}}  |u^m_{k_m}|^{n\gamma} \,
    \exp\bigg[- \frac1 2 {\rm Var} \bigg(\sum^n_{j=1} \l u_j, X(s_j,t_j)\r
    \bigg)\bigg]  \, d\overline{u} \\
    &\qquad = \frac{(2 \pi)^{(nd-1)/2}} {\big[{\rm
detCov}\big(X_0(s_1,t_1), \ldots, X_0(s_n,t_n)\big)\big]^{d/2}}\,
\int_\R
\Big(\frac v {\sigma_{m}}\Big)^{n \ga} \, e^{- \frac {v^2} 2}\, dv\\
& \qquad \le \frac{c_{_{3,11}}^n\, (n!)^\gamma} {\big[{\rm detCov}
\big( X_0(s_1,t_1), \ldots, X_0(s_n,t_n)\big)\big]^{d/2}}\, \frac1
{\sigma_{m}^{n\gamma}},
\end{split}
\end{equation}
where $\sigma_{m}^2$ is the conditional variance of
$X_{k_m}(s_m,t_m)$ given $X_i(s_j,t_j)$ ($i\ne k_m$ or $i=k_m$ but
$j\ne m$), and the last inequality follows from Stirling's formula.

Combining (\ref{Eq:M1}) and (\ref{Eq:M4}) we obtain
\begin{equation}\label{Eq:M5}
\begin{split}
&\mathcal{M}  \le \frac{c_{_{3,11}}^n\, (n!)^\gamma} {\big[{\rm
detCov} \big( X_0(s_1,t_1), \ldots, X_0(s_n,t_n)\big)\big]^{d/2}}\,
\prod_{m=1}^n \frac1 {\sigma_{m}^{\gamma}}.
\end{split}
\end{equation}

The second product in (\ref{Eq:M5}) is a ``perturbation'' factor and
will be shown to be small when integrated. For this purpose, we use
again the independence of the coordinate processes of $X$,
(\ref{Eq:XLND}) and (\ref{Eq:LND}) to derive
\begin{equation}\label{Eq:Sm1}
\begin{split}
\sigma_{m}^2 &= {\rm
Var}\Big(X_{k_m}(s_m,t_m)\Big|X_{k_m}(s_j,t_j),\
j\ne m\Big)\\
&\ge  {\rm Var}\Big(B^{\a_1}_{k_m}(s_m)\Big|B^{\a_1}_{k_m}(s_j),\
j\ne m\Big)+ {\rm
Var}\Big(B^{\a_2}_{k_m}(t_m)\Big|B^{\a_2}_{k_m}(t_j),\ j\ne
m\Big)\\
&\ge c_{_{3,12}}^2 \, \Big( \min \big\{ |s_m - s_j|^{2 \a_1}:\ j\ne
m\big\}+  \min \big\{ |t_m - t_j|^{2 \a_2}:\ j\ne m\big\}\Big).
\end{split}
\end{equation}

As in the proof of Eq. (\ref{Eq:momentEst1}), we will prove Eq.
(\ref{Eq:momentEst2}) by cases.

If $\frac{N_1}{\a_1}>d$, then we take
$\ga\in\big(0,\,\frac12(\frac{N_1}{\a_1}-d)\big)$  so that
\begin{equation} \label{Eq:d21ga}
\a_1(d+2\gamma)<N_1.
\end{equation}

For any $n$ points $(s_1,t_1), \ldots, (s_n,t_n) \in D\backslash
\{0\}$, we define a permutation $\pi_s$ of $\{1, 2, \ldots, n\}$
such that
\begin{equation}\label{Eq:perm1}
\begin{split}
&\quad |s_{\pi_s(1)}|=\min\{|s_i|,\,i=1,\ldots,n\}, \nonumber\\
&\quad
|s_{\pi_s(j)}-s_{\pi_s(j-1)}|=\min\left\{|s_i-s_{\pi_s(j-1)}|,\,
i\in\{1,\ldots,n\}\backslash\{\pi_s(1),\ldots,\pi_s(j-1)\}\right\}.\nonumber
\end{split}
\end{equation}
Then, by (\ref{Eq:Sm1}), we have
\begin{equation}\label{Eq:Prods}
\begin{split}
\prod_{m=1}^n \frac1 {\sigma_{m}^{\gamma}} &\le \prod_{m=1}^n
\frac{1} {c_{_{3,12}}\, \big[|s_{\pi_s(m)} - s_{\pi_s(m-1)}|^{\a_1}
\wedge |s_{\pi_s(m+1)} -
s_{\pi_s(m)}|^{\a_1} \big]^{\gamma}}\\
&\le c_{_{3,12}}^{-n} \prod_{m=1}^n \frac1 {|s_{\pi_s(m)} -
s_{\pi_s(m-1)}|^{2 \a_1 \gamma}}\\
&\le c_{_{3,12}}^{-n} \, \frac 1 {\big[{\rm detCov}
\big(B^{\a_1}_0(s_1),\ldots, B^{\a_1}_0(s_n) \big)\big]^{\ga}}.
\end{split}
\end{equation}

It follows from (\ref{Eq:M5}), (\ref{Eq:Det_XB}) and
(\ref{Eq:Prods}) that
\begin{equation}\label{Eq:M61}
\begin{split}
\mathcal{M} & \le  \frac {c_{_{3,11}}^n\, (n!)^\gamma} {\big[{\rm
detCov} \big(B^{\a_1}_0(s_1),\ldots, B^{\a_1}_0(s_n)
\big)\big]^{d/2}}\, \prod_{m=1}^n \frac1 {\sigma_{m}^{\gamma}}\\
&\le c_{_{3,13}}^n\, (n!)^\gamma\, \frac {1} {\big[{\rm detCov}
\big(B^{\a_1}_0(s_1),\ldots, B^{\a_1}_0(s_n)
\big)\big]^{(d+2\ga)/2}}\\
&\le c_{_{3,14}}^n\, (n!)^\gamma\, \prod_{j=1}^n \frac {1} { \min
\{|s_j - s_i|^{\a_1 (d+ 2\ga)}, 0 \le i \le j-1\}}.
\end{split}
\end{equation}

Therefore, by (\ref{Eq:d21ga}) and Lemma \ref{Lem:Xiao97}, we have
\begin{equation}\label{Eq:M91}
\begin{split}
\int_{D^n} \mathcal{M}( \overline k, \overline{\t}, \gamma) \,
d\overline{\t} &\le  c_{_{3,14}}^n\, (n!)^\gamma\, \int_{D^n}
\prod_{j=1}^n \frac {1} { \min
\{|s_j - s_i|^{\a_1 (d+ 2\ga)}, 0 \le i \le j-1\}}\, d\overline{\t}\\
&\le c_{_{3,15}}^n\,
(n!)^{\frac{\a_1(d+2\gamma)}{N_1}+  \gamma }\, r^{n\big(N_1-\a_1(d +
2\gamma)\big)}\times r^{nN_2}\\
&= c_{_{3,15}}^n\, (n!)^{\eta_1+ \big(\frac{2\a_1}{N_1}+1\big)
\gamma }\, r^{n(\beta_1- 2\a_1\gamma)}.
\end{split}
\end{equation}

We combine (\ref{Eq:LTmoment3}) and (\ref{Eq:M91}) to obtain
\begin{equation} \label{Eq:LTmoment4}
\begin{split}
\E \Big[ \big(L(x, D) - L(y, D)\big)^n\Big]\le c_{_{3,16}}^n\,
(n!)^{\eta_1+ \big(\frac{2\a_1}{N_1}+1\big) \gamma
}\,|x-y|^{n\ga}\, r^{n(\beta_1-2\a_1\gamma)}.
\end{split}
\end{equation}
By choosing the constant $\kappa_1 \ge \max\{\frac{2\a_1}{N_1}+1,
2\a_1\}$, we prove Eq. (\ref{Eq:momentEst2}) for the case
$\frac{N_1}{\a_1}>d$ [i.e., $\tau=1$].

Now we prove Eq. (\ref{Eq:momentEst2}) for the case of
$\frac{N_1}{\a_1}\le d<\frac{N_1}{\a_1}+\frac{N_2}{\a_2}$. Inspired
by Lemma 3.4 in Ayache, Wu and Xiao (2008), we choose
$$
\ga\in\bigg(0,\, \frac1
4\Big(\frac{N_1}{\a_1}+\frac{N_2}{\a_2}-d\Big)\bigg),
$$
\begin{equation}\label{Eq:Ho_de}
\delta=\frac12\min\left\{1,\,\a_1\big(\frac{N_1}{\a_1}+
\frac{N_2}{\a_2}-d\big),\,\a_1 \ga \right\}
\end{equation}
and set
\begin{equation}\label{Eq:Ho_p}
\frac1{p_1}=\frac{N_1-\delta}{\a_1d},\qquad\frac1{p_2}=1-\frac1{p_1}.
\end{equation}
Clearly, we have
\begin{equation}\label{Eq:Ho00}
p_1>1,\ \,p_2>1,\ \, \ \ \frac1{p_1}+\frac1{p_2}=1
\end{equation}
and
\begin{equation}\label{Eq:Ho10}
\begin{split}
&\frac{\a_1 d}{p_1}=N_1-\delta<N_1,\\
&\frac{\a_2
d}{p_2}=\a_2\Big(d-\frac{N_1}{\a_1}+\frac{\delta}{\a_1}\Big)<N_2,
\end{split}
\end{equation}
where the last inequality follows from the fact that $\delta\le
\frac{\a_1}2\big(\frac{N_1}{\a_1}+\frac{N_2}{\a_2}-d\big).$ By a
simple computation, we also have
\begin{equation}\label{Eq:Ho11}
\frac{\a_1 d}{p_1}+\frac{\a_2
d}{p_2}=N_1+\a_2d-\frac{\a_2}{\a_1}N_1+\Big(\frac{\a_2}
{\a_1}-1\Big)\delta =N-\beta_2+\Big(\frac{\a_2}{\a_1}-1\Big)\delta
\end{equation}
and
\begin{equation}\label{Eq:Ho12}
\frac{\a_1 d}{N_1p_1}+\frac{\a_2 d}{N_2p_2}=1+\frac{\a_2
d}{N_2}-\frac{\a_2 N_1}{\a_1 N_2}+ \Big(\frac{\a_2
}{\a_1N_2}-\frac1{N_1}\Big)\delta=\eta_2+\Big(\frac{\a_2
}{\a_1N_2}-\frac1{N_1}\Big)\delta.
\end{equation}
Furthermore, from the way we define $\ga,\, \delta$ and $p_2$, we
know
\begin{equation}\label{Eq:Ho_ga}
\frac{\a_2 d}{p_2}+2\a_2\ga<N_2.
\end{equation}

For any $n$ points $(s_1,t_1), \ldots, (s_n,t_n) \in D\backslash
\{0\}$, we define a permutation $\pi_t$ of $\{1, 2, \ldots, n\}$
such that
\begin{equation}\label{Eq:perm2}
\begin{split}
&\quad |t_{\pi_t(1)}|=\min\{|t_i|,\,i=1,\ldots,n\}, \nonumber\\
&\quad
|t_{\pi_t(j)}-t_{\pi_t(j-1)}|=\min\left\{|t_i-t_{\pi_t(j-1)}|,\,
i\in\{1,\ldots,n\}\backslash\{\pi_t(1),\ldots,\pi_t(j-1)\}\right\}.\nonumber
\end{split}
\end{equation}
Then, by (\ref{Eq:Sm1}), we have
\begin{equation}\label{Eq:Prodt}
\begin{split}
\prod_{m=1}^n \frac1 {\sigma_{m}^{\gamma}} &\le \prod_{m=1}^n
\frac{1} {c_{_{3,12}}\, \big[|t_{\pi_t(m)} - t_{\pi_t(m-1)}|^{\a_2}
\wedge |t_{\pi_t(m+1)} -
t_{\pi_t(m)}|^{\a_2} \big]^{\gamma}}\\
&\le c_{_{3,12}}^{-n}\, \frac 1 {\big[{\rm
detCov}\left(B^{\a_2}_0(t_1),\ldots,,
B^{\a_2}_0(t_n)\right)\big]^{\ga}}.
\end{split}
\end{equation}

Recall from (\ref{Eq:Det_XB}) that
\begin{equation}\label{Eq:Det_XB-2}
\begin{split}
{\rm detCov}&\left(X_0(s_1,t_1),\ldots,,X_0(s_n,t_n)\right)\\
&\ge {\rm
detCov}\left(B^{\a_1}_0(s_1),\ldots,,B^{\a_1}_0(s_n)\right) +{\rm
detCov}\left(B^{\a_2}_0(t_1),\ldots,,B^{\a_2}_0(t_n)\right).
\end{split}
\end{equation}
Hence,
\begin{equation}\label{Eq:Det_XB-3}
\begin{split}
\Big[{\rm detCov}\left(X_0(s_1,t_1), \ldots,
X_0(s_n,t_n)\right)\Big]^{-\frac 1 2} & \le  \Big[{\rm
detCov}\left(B^{\a_1}_0(s_1),
 \ldots, B^{\a_1}_0(s_n)\right)\Big]^{-\frac 1 {2p_1}}\\
&\quad \times
 \Big[{\rm detCov}\left(B^{\a_2}_0(t_1),
 \ldots, B^{\a_2}_0(t_n)\right)\Big]^{-\frac 1 {2p_2}}.
 \end{split}
\end{equation}
It follows from (\ref{Eq:M5}), (\ref{Eq:Prodt}) and
(\ref{Eq:Det_XB-3})
 that
\begin{equation}\label{Eq:M6}
\begin{split}
\mathcal{M} & \le \, \, \frac {c_{_{3,11}}^n\, (n!)^\gamma}
{\big[{\rm detCov} \big(X_0(s_1, t_1),\ldots, X_0(s_n, t_n)
\big)\big]^{d/2} } \, \prod_{m=1}^n \frac1 {\sigma_{m}^{\gamma}}\\
&\le \frac {c_{_{3,17}}^n\, (n!)^\gamma} {\big[{\rm detCov}
\big(B^{\a_1}_0(s_1),\ldots, B^{\a_1}_0(s_n)
\big)\big]^{\frac{d}{2p_1}}\,\big[{\rm
detCov}\left(B^{\a_2}_0(t_1),\ldots,,
B^{\a_2}_0(t_n)\right)\big]^{\frac{d} {2p_2} +  \gamma }}.
\end{split}
\end{equation}

Combining (\ref{Eq:M6}), (\ref{Eq:Ho10}), (\ref{Eq:Ho_ga}) and
Lemma \ref{Lem:Xiao97}, we obtain
\begin{equation}\label{Eq:M9}
\begin{split}
\int_{D^n} \mathcal{M}( \overline k, \overline{\t}, \gamma) \,
d\overline{\t} &\le c_{_{3,17}}^n\,(n!)^\gamma\,
\int_{O_{N_1}^n(u_1,r)} \frac{d \overline{s}} {\big[{\rm detCov}
\big(B^{\a_1}_0(s_1),\ldots, B^{\a_1}_0(s_n)
\big)\big]^{\frac{d}{2p_1}}} \\
&\qquad \qquad \times \int_{O_{N_2}^n(u_2,r)} \frac{d \overline{t}}
{\big[{\rm detCov} \big(B^{\a_2}_0(t_1),\ldots, B^{\a_2}_0(t_n)
\big)\big]^{\frac{d}{2p_2} + \gamma}}\\
&\le c_{_{3,18}}^n\, (n!)^{\sum_{\ell=1}^2\frac {\a_\ell\, d}{N_\ell
p_\ell}+\big(1+\frac{2\a_2}{N_2}\big)\gamma} r^{n
\left(N-\sum_{\ell=1}^2\frac{\a_\ell\,d
}{p_\ell}-\a_2\gamma\right)}\\
&\le c_{_{3,18}}^n\,(n!)^{\eta_2 + \kappa_1\gamma}\,  r^{n (\beta_2-
\kappa_1 \gamma)}.
\end{split}
\end{equation}
In the above, the constant $\kappa_1 >0$ is chosen appropriately by
taking into account (\ref{Eq:Ho11}), (\ref{Eq:Ho12}) and
(\ref{Eq:LTmoment4}). The value of $\kappa_1$ depends on $\a_1,
\a_2, N_1, N_2$ and $d$ only.

We combine (\ref{Eq:LTmoment3}) and (\ref{Eq:M9}) to obtain
\begin{equation} \label{Eq:LTmoment41}
\E \Big[ \big(L(x, D) - L(y, D)\big)^n\Big] \le
c_{_{3,9}}^n\,(n!)^{\eta_2 + \kappa_1
 \gamma}\,|x-y|^{n\ga}  r^{n
(\beta_2- \kappa_1 \gamma)}.
\end{equation}
This proves (\ref{Eq:momentEst2}) for the case of
$\frac{N_1}{\a_1}\le d<\frac{N_1}{\a_1}+\frac{N_2}{\a_2}$ [i.e.,
$\tau=2$]. The proof of Lemma \ref{Lem:MomEst2} is complete.
\end{proof}

Now we are ready to prove Theorem \ref{Thm:JCont}.

\begin{proof}{\bf of Theorem \ref{Thm:JCont}} The proof of the joint
continuity of the local time of $X$ is similar to that of Theorem 3.1
in Ayache, Wu and Xiao (2008) [see also the proof of Theorem 8.2 in
Xiao (2009)].  Hence we only give a sketch of it.

It suffices to show that for any fixed $u= (u_1, u_2) \in \R^N$ and
$R > 0$, the local time $L(x, (s, t)):= L(x, [u_1, u_1+s] \times
[u_2, u_2+t])$ has a version which is continuous in $(x, s, t) \in
\R^d \times [0,\,R]^N$ almost surely. For simplicity of notation, we
assume $u = (0, 0)$. Observe that for all $x, y \in \R^d$,
$(s,t),\, (v,w) \in [0,\,R]^N$ and all even integers $n \ge 1$, we
have
\begin{equation}\label{Eq:Last}
\begin{split}
\E\Big[&\big(L(x,[0,\,s]\times [0,\,t]) - L(y, [0,\,v]
\times [0,\,w])\big)^n \Big]\\
&\le 2^{n-1} \bigg\{\E\Big[\big(L(x, [0,\,s]\times [0,\,t]) -
L(x, [0,\,v]\times [0,\,w])\big)^n \Big]\\
&\qquad+ \E\Big[\big(L(x, [0,\,v]\times [0,\,w]) - L(y,
[0,\,v]\times [0,\,w])\big)^n\Big]
\bigg\}.\\
\end{split}
\end{equation}
Since $L(x, \cdot)$ is a finite Borel measure, the difference $L(x,
[0,\,s]\times [0,\,t]) - L(x, [0,\,v]\times [0,\,w])$ can be bounded
by a sum of finite number of terms of the form $L(x, D_j)$, where
each $D_j$ is a closed subset of $[0,\,R]^N$ of the form
$O_{N_1}(\cdot,r)\times O_{N_2}(\cdot,r)$ with the radius $r\le
\frac12|(s,t)-(v,w)|:=\frac12\sqrt{|s-v|^2+|t-w|^2}.$ We can use
(\ref{Eq:momentEst1}) to bound the first term in (\ref{Eq:Last}). On
the other hand, the second term in (\ref{Eq:Last}) can be dealt with
by using (\ref{Eq:momentEst2}). Consequently, there exist
some constants $\gamma \in (0, 1)$ and $n_0$ such that for all $x, y
\in \R^d$,  $(s,t),\, (v,w) \in O_{N_1,N_2}(0,R)\cap\R^N_+$ and all
even integers $n \ge n_0$,
\begin{equation}\label{Eq:Last2}
\begin{split}
\E\Big[&\big(L(x,[0,\,s]\times [0,\,t]) - L(y, [0,\,v]\times [0,\,w])\big)^n \Big]\\
& \le c_{_{3,19}}^n\,\big(|x-y| + |(s,t)-(v,w)|\big)^{n\ga}.
\end{split}
\end{equation}
It follows from (\ref{Eq:Last2}) and  the multiparameter version of
Kolmogorov's continuity theorem [cf. Khoshnevisan (2002)] that
there exists a modification of the local times of $X$, still denoted
by $L(x, (s,t))$, such that it is continuous for
$x\in\R^d,\,(s,t)\in [0,\, R]^N$. This finishes the proof of Theorem
\ref{Thm:JCont}.
\end{proof}


\section{Exponential integrability and H\"older conditions
for the intersection local times} \label{Sec:H}

In this section, we investigate the exponential integrability and
asymptotic behavior of the local time $L(x, \cdot)$ of $X$. As
applications of the later result, we obtain a lower bound for
the exact Hausdorff measure of the set $M_2$ of the intersection
times of $B^{\a_1}$ and $B^{\a_2}$.

The following two technical lemmas will play essential r\^oles in
our derivation.

\begin{lemma} \label{Lem:taubounds}
Under the conditions of
Theorem \ref{Thm:JCont}, there exist positive and finite constants $\eps \in (0,
1/e)$, $c_{_{4,1}}$ and $c_{_{4,2}}$, depending on $\a_1,\, \a_2$,
$N_1,\,N_2$, and $d$ only, such that the following hold:
\begin{itemize}
\item[(i)]\ For all $(a_1,a_2) \in \R^N$ and
$D = O_{N_1,N_2}\big((a_1,a_2),r\big) $ with radius $r \in (0,
\eps)$, $x \in \R^d$ and all integers $n \ge 1$,
\begin{equation}\label{Eq:m11x}
 \E \Big[ L\big(x+X(a_1,a_2), D \big)^n\Big] \le  \left\{\begin{array}{ll}
c_{_{4,1}}^n\, (n!)^{\eta_\tau}\,
r^{n\,\beta_\tau} \quad &\hbox{ if } \frac{N_1}{\a_1} \ne d,\\
c_{_{4,1}}^n\, n!\, r^{n\,N_2} \, \prod_{j=1}^n\log \bigg(e+ \frac
{j^{\left(\frac{\a_2} {N_2} - \frac {\a_1} {N_1}\right)^+}}
{r^{\a_2-\a_1}} \bigg)\quad &\hbox{ if } \frac{N_1}{\a_1}=d.
\end{array}\right.
\end{equation}

\item[(ii)]\ For all $(a_1,a_2) \in \R^N$ and
$D = O_{N_1,N_2}\big((a_1,a_2),r\big)$ with radius $r >0$, $x,\, y
\in \R^d$ with $|x - y|\le 1$,  all even integers $n \ge 1$ and all
$\gamma \in (0, 1)$ small,
\begin{equation}\label{Eq:m20y}
\begin{split}
\E \Big[\big(L(x & +X(a_1,a_2), D) - L(y+ X(a_1,a_2), D)\big)^n\Big] \\
& \le c_{_{4,2}}^n\, (n!)^{\eta_\tau + \kappa_1 \gamma}\, |x-y|^{n
\gamma}\, r^{n (\beta_\tau - \kappa_1 \gamma)}.
\end{split}
\end{equation}
In the above, $\kappa_1>0$ is the same constant as in Lemma \ref{Lem:MomEst2}.
\end{itemize}
\end{lemma}

\begin{proof}\ For any fixed $(a_1,a_2) \in \R^N$, we define the
Gaussian random field $Y =\{Y(s,t),$ $ (s,t) \in \R^N\}$ with values
in $\R^d$ by $Y(s,t) = X(s,t) - X(a_1,a_2)$. It follows from
(\ref{Eq:occupation}) that if  $X$ has a local time $L(x,S)$ on any
Borel set $S$, then $Y$ also has a local time $\tilde{L}(x, S)$ on
$S$ and, moreover, $L(x + X(a_1,a_2),S ) = \tilde{L}(x, S)$. Since
$X$ has stationary increments, both Lemma \ref{Lem:MomEst} and
Lemma \ref{Lem:MomEst2} hold
for the Gaussian field $Y$. This proves (\ref{Eq:m11x}) and
(\ref{Eq:m20y}).
\end{proof}

The following lemma is a consequence of Lemma \ref{Lem:taubounds}
and Chebyshev's inequality.

\begin{lemma}\label{Lem:probbounds}
Assume  the conditions of
Theorem \ref{Thm:JCont} hold. For any $b >0$, there exist
positive and finite constants $\eps \in (0, 1/e)$,
$c_{_{4,3}}, \, c_{_{4,4}},\, c_{_{4,5}}$, (depending on
$\a_1,\, \a_2,\,N_1,\,N_2$ and $d$ only), such that for
all $(a_1, a_2)\in \R^N$, $D=O_{N_1,N_2}\big((a_1,a_2),
r\big)$ with $r \in (0, \,\eps)$, $x\in \R^d$ and $u >1 $
large enough, the following inequalities hold:
\begin{itemize}
\item[(i).] If $N_1 \ne \a_1 d$, then
\begin{equation}\label{Eq:ProbB1}
\P\bigg\{L\big(x+X(a_1,a_2), \, D\big) \ge \,c_{_{4,3}}\,
r^{\beta_\tau} \, u^{\eta_\tau} \bigg\} \le \exp \big({-b\, u}
\big).
\end{equation}
\item[(ii).] If $N_1 = \a_1 d$, then
\begin{equation}\label{Eq:ProbB2}
\P\bigg\{L\big(x+X(a_1,a_2), \, D\big) \ge \,c_{_{4,4}}\, r^{N_2} \,
u\, \log\Big(e+\frac {u^{\left(\frac{\a_2} {N_2} - \frac {\a_1}
{N_1}\right)^+}} {r^{\a_2-\a_1}} \Big) \bigg\} \le \exp \big({-b\,
u} \big).
\end{equation}
\item[(iii).] For $x,\, y \in \R^d$ with $|x - y|\le 1$  and
$\gamma>0$ small,
\begin{equation}\label{Eq:ProbB3}
\begin{split}
\P\bigg\{\Big| L\big(x &+X(a_1,a_2),\, D \big)-
L\big(y+ X(a_1,a_2),\, D\big)\Big|\\
& \qquad \ \ge c_{_{4,5}}\, |x-y|^{\gamma}\, r^{\beta_\tau- \kappa
\gamma}\, u^{\eta_\tau+\kappa \gamma} \bigg\} \le \exp\big({-b\, u}
\big).
\end{split}
\end{equation}
\end{itemize}
\end{lemma}

\begin{proof}\ The proofs of Parts (i) and (iii) based on Lemma
\ref{Lem:taubounds} and
Chebyshev's inequality are standard, hence omitted. In the following
we prove (ii). Define the random variable $\Lambda =
L\big(x+X(a_1,a_2), D\big)/r^{N_2}$. For $u > 0$ large, let $n =
\lfloor u \rfloor$, the largest positive integer no bigger than $u$.
We apply Chebyshev's inequality and Lemma \ref{Lem:taubounds} to obtain
\begin{equation}\label{Eq:ProbB4}
\P\Bigg\{ \Lambda \ge \,c \, u \log\bigg(e+\frac
{u^{\left(\frac{\a_2} {N_2} - \frac {\a_1} {N_1}\right)^+}}
{r^{\a_2-\a_1}} \bigg) \Bigg\} \le \bigg(\frac{ c_{_{4,1}}} {e
c}\bigg)^n  \, \frac{ \prod_{j=1}^n\log \Big(e+\frac
{j^{\left(\frac{\a_2} { N_2} - \frac {\a_1} {N_1}\right)^+}}
{r^{\a_2- \a_1}} \Big)} {  \log^n\Big(e+\frac {n^{\left(\frac{\a_2}
{ N_2} - \frac {\a_1} {N_1}\right)^+}} {r^{\a_2- \a_1}} \Big)} \le
\bigg(\frac{ c_{_{4,1}}} {e c}\bigg)^n,
\end{equation}
where $c>0$ is a constant whose value will be determined later, and
where we have used the fact that for $j\in\{1,\,2,\ldots,n\},$
\[
\log \bigg(e+\frac {j^{\left(\frac{\a_2} { N_2} - \frac {\a_1}
{N_1}\right)^+}} {r^{\a_2- \a_1}} \bigg)\le \log \bigg(e+\frac
{n^{\left(\frac{\a_2} { N_2} - \frac {\a_1} {N_1}\right)^+}}
{r^{\a_2- \a_1}} \bigg).
\]
By taking $c = c_{_{4,4}}$ large so that $\log \big( c_{_{4,1}}/(e
c)\big) \le -b$, we obtain  (\ref{Eq:ProbB2}).
\end{proof}

The following result about the exponential integrability of $L(x,
D)$ is a direct consequence of Lemma \ref{Lem:probbounds}. We
omit its proof.

\begin{theorem}\label{Thm:Exp}
Assume that the conditions of Theorem \ref{Thm:JCont} hold
and let $D_1:=O_{N_1,N_2}(0,1)$. Then there exists a constant
$\delta > 0$, depending on $\a_1$, $\a_2$, $N_1, N_2$ and $d$ only,
such that the following hold:
\begin{itemize}
\item[(i).] If $N_1 \ne \a_1 d$, then for every $x \in \R^{d}$
\begin{equation}\label{Eq:m20}
\E \Big(e^{\delta L(x, D_1)^{\eta_\tau}}\Big) \le 1,
\end{equation}
where  $\eta_\tau$ is the constant given in (\ref{Def:gamma}).
\item[(ii).] If $N_1 = \a_1 d$, then
\begin{equation}\label{Eq:m20-b}
\E \Big(e^{\delta \psi_3(L(x, D_1))}\Big) \le 1,
\end{equation}
where $\psi_3(y) = y/\log(e+ y)$ for all $y > 0$.
\end{itemize}
\end{theorem}

Now we study the local H\"older condition of the intersection local
time $L^{\a_1,\a_2}(\cdot)$ and its connection to fractal properties
of the set of intersection times $M_2$ of $B^{\a_1}$ and $B^{\a_2}$.

Since $M_2$ is the zero-set of $X$, namely, $M_2 = X^{-1}(0)$, and
the Gaussian random field $X$ satisfies the conditions in
Xiao (2009). It follows from  Theorem 7.1 in Xiao (2009) that
\begin{equation}\label{Eq:dimM20}
\dim M_2 = \dimp M_2 = \beta_\tau
\end{equation}
with positive probability. In the above, $\dim$ and $\dimp$ denote
Hausdorff and packing dimension, respectively; see Falconer (1990)
for more information. In Corollary \ref{Cor:dimM2} below, we will
show that (\ref{Eq:dimM20}) holds with probability 1.

It is an interesting problem to determine the exact Hausdorff and
packing measure functions for $M_2$. For this purpose, the limsup
and liminf type laws of the iterated logarithm need to be
established, respectively for the intersection local time
$L^{\a_1,\a_2}(\cdot)$.

In the following, we consider the limsup laws of the iterated logarithm
for the local time $L(x,\cdot)$ of $X$.
By applying Lemma \ref{Lem:probbounds} [with $(a_1,a_2) =(0,0)$] and the
Borel-Cantelli lemma, one can easily derive the following result: There
exists a positive constant $c_{_{4,6}}$ such that for every $x \in \R^d$
and $(s,t) \in \R^N$,
\begin{equation}\label{Eq:Holder0}
\limsup_{r \to 0} \frac{L\big(x, O_{N_1,N_2}((s,t), r)\big)}
{\varphi_1(r)} \le c_{_{4,6}}, \qquad \hbox{ a.s.},
\end{equation}
where
\begin{equation}\label{Eq:phi1}
\varphi_1(r) = \left\{\begin{array}{ll} r^{\beta_\tau}\, \big(\log
\log (1/r)\big)^{\eta_\tau}
\qquad &\hbox{ if } N_1 \ne \a_1 d,\\
r^{N_2} \big(\log \log (1/r)\big)\,\log \Big(e+ \frac{(\log \log
(1/r))^{\left(\frac{\a_2} {N_2}
 - \frac {\a_1} {N_1}\right)^+}} {r^{\a_2 - \a_1}}\Big)
&\hbox{ if } N_1 = \a_1 d.
\end{array}
\right.
\end{equation}
It is worthwhile to compare (\ref{Eq:Holder0}) with the
corresponding results for $(N,d)$ fractional Brownian motion of
index $\a$ in Xiao (1997) and the $(N, d)$ fractional Brownian
sheets with index $(\a, \ldots, \a)\in (0, 1)^N$ in Ayache, Wu and
Xiao (2008). In the former case, $X$ is isotropic and its local
time $L(x, \cdot)$ satisfies
\begin{equation}\label{Eq:Holder-fbm}
\limsup_{r \to 0} \frac{L\big(x, O_{N}(t, r)\big)} {r^{N-\a d}
(\log\log 1/r)^{\a d/N}} \le c_{_{4,7}}, \qquad \hbox{ a.s.},
\end{equation}
while the local time of the $(N, d)$ fractional Brownian sheet with
index $(\a, \ldots, \a) \in (0, 1)^N$ satisfies
\begin{equation}\label{Eq:Holder-sh}
\limsup_{r \to 0} \frac{L\big(x, O_{N}(t, r)\big)} {r^{N- \a d}
(\log\log 1/r)^{\a d}} \le c_{_{4,8}}, \qquad \hbox{ a.s.}
\end{equation}
Note that, the anisotropy of the fractional Brownian sheet only
increases the power of the correction factor $\log\log 1/r$. For the
Gaussian random field $X$ defined by (\ref{Def:X}) with $N_1 = \a_1
d$, (\ref{Eq:Holder0}) suggests that the asymptotic properties of
the local times of $X$ may be significantly different from those in
(\ref{Eq:Holder-fbm}) and (\ref{Eq:Holder-sh}). In fact, when $N_1 = \a_1
d$ and as $r\downarrow 0$, we have
\[\varphi_1(r)\sim\left\{\begin{array}{ll}
r^{N_2} \log \log (1/r)\,\log \log \log (1/r)
&\hbox{ if } \a_2=\a_1,\,\,\frac{\a_2}{N_2}>\frac{\a_1}{N_1}, \\
r^{N_2}\log \log (1/r) &\hbox{ if } \a_2=\a_1,\,\,\frac{\a_2}{N_2}
\le\frac{\a_1}{N_1},\\
r^{N_2}\log(1/r)\log\log(1/r)  &\hbox{ if } \a_2>\a_1.
\end{array}
\right.\]
However, in this later case, it is unclear to us whether
the logarithmic correction factor in (\ref{Eq:phi1}) is sharp. It would be
interesting to study this problem and establish sharp laws of the
iterated logarithm for the local times of $X$. For such a result for
the local times of a one-parameter fractional Brownian motion, see
Baraka and Mountford (2008).

As a consequence of (\ref{Eq:Holder0}) we have for the intersection
local time of $B^{\a_1}$ and $B^{\a_2}$ that, for every $(s,t) \in
\R^N$,
\begin{equation}\label{Eq:Holder01}
\limsup_{r \to 0} \frac{L^{\a_1,\a_2}\big(O_{N_1,N_2}((s,t),
r)\big)} {\varphi_1(r)} \le c_{_{4,6}}, \qquad \hbox{ a.s.}
\end{equation}
It follows from Fubini's theorem that, with probability one,
(\ref{Eq:Holder01}) holds for $\lambda_N$-almost all $(s,t) \in
\R^N$. Now we prove a stronger version of this result, which is
useful in determining the exact Hausdorff measure of $M_2$.

\begin{theorem}\label{Th:Holder3}
Assume that $
\frac{N_1}{\a_1}+\frac{N_2}{\a_2}>d$. Let $\tau \in \{1, \,2\}$ be
the integer defined in (\ref{Def:tau}) and let $D =
O_{N_1,N_2}(0,R)$ be fixed. Let $L^{\a_1,\a_2}(\cdot)$ be the
intersection local time of $B^{\a_1}$ and $B^{\a_2}$, which is a
random measure supported on the set $M_2$. Then there exists a
positive and finite constant $c_{_{4,9}}$ such that with probability
$1$,
\begin{equation}\label{Eq:Holder3}
\limsup_{r \to 0} \frac{L^{\a_1,\a_2}\big(O_{N_1,N_2}((s,t),
r)\big)} {\varphi_1(r)} \le c_{_{4,9}}
\end{equation}
holds for $L^{\a_1,\a_2}(\cdot)$-almost all $(s,t)\in D$, where
$\varphi_1(r)$ is defined in (\ref{Eq:phi1}).
\end{theorem}

\begin{proof}\ Again we work on the random field $X$ defined by
(\ref{Def:X}). For every integer $k > 0$, we consider the random
measure $L_k(x, \bullet)$ on the Borel subsets $C$ of
$O_{N_1,N_2}(0,R)$ defined by
\begin{equation}\label{Eq:mun}
\begin{split}
L_k(x,C) &= \int_C (2\pi k)^{d/2} \exp\bigg( - \frac{k\,
|X(s,t)-x|^2} 2\bigg)\,ds dt \\
&= \int_C \int_{\R^d} \exp \bigg(-\frac{|\xi|^2} {2k}  + i\langle
\xi, X(s,t)-x\rangle \bigg) d\xi \, dsdt.
\end{split}
\end{equation}
Then, by the occupation density formula (\ref{Eq:occupation}) and
the continuity of the function $y \mapsto L(y, C)$, one can verify
that almost surely $L_k(x, C) \to L(x, C)$ as $k \to \infty$ for
every Borel set $C \subseteq O_{N_1,N_2}(0,R)$.

For every integer $m \ge 1$, denote $f_m(s,t) = L\big(x,\,
O_{N_1,N_2}((s,t),2^{-m})\big)$. From the proof of Theorem \ref{Thm:JCont} we
can see that almost surely the functions $f_m(s,t)$ are continuous
and bounded. Hence we have almost surely, for all integers $m,\, n
\ge 1$,
\begin{equation}\label{Eq:Approx}
\int_{O_{N_1,N_2}(0,R)} \left[f_m(s,t)\right]^n\, L(x,ds\,dt) =
\lim_{k \to \infty} \int_{O_{N_1,N_2}(0,R)}
\left[f_m(s,t)\right]^n\, L_k(x,ds\,dt).
\end{equation}
It follows from (\ref{Eq:Approx}), (\ref{Eq:mun}) and the proof of
Proposition 3.1 of Pitt (1978) that, for every positive integer $
n\ge 1$,
\begin{equation}\label{Eq:Holder31}
\begin{split}
&\E\int_{O_{N_1,N_2}(0,R)} \left[f_m(s,t)\right]^n\, L(x,dsdt)\\
&\quad=\left(\frac 1 {2\pi}\right)^{(n+1)d} \int_{O_{N_1,N_2}(0,R)}
\int_{O_{N_1,N_2}((s,t),2^{-m})^n} \int_{\R^{(n+1)d}}
\exp\bigg(-i\sum_{j=1}^{n+1}\l x,u^j\r\bigg)
\\ &\qquad \qquad\qquad\qquad \times
\E\exp\bigg(i\sum_{j=1}^{n+1}\l u^j, X(s_j,t_j)\r\bigg)\,
d{\overline u}d{\overline \t},
\end{split}
\end{equation}
where ${\overline u}=(u^1,\ldots,u^{n+1})\in\R^{(n+1)d}$ and
${\overline \t}=(s,\,t,\,s_1,\,t_1,\ldots,s_n,\,t_n)$. Similar to
the proof of (\ref{Eq:momentEst1}) we have that the right hand side
of (\ref{Eq:Holder31}) is at most
\begin{equation}\label{Eq:Holder32}
\begin{split}
&\int_{O_{N_1,N_2}(0,R)}\int_{O_{N_1,N_2}((s,t),2^{-m})^n}
\frac{c_{_{4,10}}^n\, d{\overline \t}}{\Big[{\rm
detCov}\big(X_0(s,t),X_0(s_1,t_1),\ldots,X_0(s_n,t_n)\big)\Big]^{d/2}}\\
& \le \left\{\begin{array}{ll}
c_{_{4,11}}^n\,(n!)^{\eta_\tau}\, 2^{-mn\beta_\tau}, \qquad &\hbox{ if } N_1\ne \a_1 d,\\
c_{_{4,11}}^n\, n!\, 2^{-nm\,N_2} \, \prod_{j=1}^n\log \Big(e+
j^{\left(\frac{\a_2} {N_2} -\frac {\a_1} {N_1}\right)^+ } 2^{(\a_2-
\a_1)m} \Big) \qquad &\hbox{ if } N_1= \a_1 d,
\end{array}\right.
\end{split}
\end{equation}
where $c_{_{4,11}}$ is a positive finite constant depending on
$\a_1,\,\a_2,N_1,\,N_2,\,d$ and $R$ only.

Let $ \rho>0$ be a constant whose value will be determined later. We
consider the random set
\[
D_m(\omega)=\left\{(s,t)\in O_{N_1,N_2}(0,R):\, f_m(s,t) \ge \rho\,
\varphi_1(2^{-m})\right\}.
\]
Denote by $\mu_\omega$ the restriction of the random measure
$L(x,\cdot)$ on $O_{N_1,N_2}(0,R)$, that is, $\mu_\omega(E)=L(x,\,
E\cap O_{N_1,N_2}(0,R))$ for every Borel set $E\subseteq \R^N$. Now we
take $n= \lfloor \log m \rfloor$. Then, by applying
(\ref{Eq:Holder32}) and by Stirling's formula, we have
\begin{equation} \label{Eq:Uptail}
\begin{split}
\E\mu_\omega(D_m)&\le
\frac{\E\int_{O_{N_1,N_2}(0,R)}\left[f_m(s,t)\right]^n\,
L(x,ds\,dt)} {[ \rho\, \varphi_1(2^{-m})]^n} \le m^{-2},
\end{split}
\end{equation}
provided $\rho >0$ is chosen large enough, say, $\rho \ge
c_{_{4,11}}\, e^2:= c_{_{4,9}}$. This implies that
\[
\E\left(\sum_{m=1}^\infty\mu_\omega(D_m)\right)<\infty.
\]
Therefore, with probability $1$ for $\mu_\omega$ almost all
$(s,t)\in O_{N_1,N_2}(0,R)$, we derive
\begin{equation}\label{Eq:Holder33}
\limsup_{m \to \infty} \frac{L(x, O_{N_1,N_2}((s,t), 2^{-m}))}
{\varphi_1(2^{-m})} \le c_{_{4,9}}.
\end{equation}

Finally, for any $r>0$ small enough, there exists an integer $m$
such that $2^{-m}\le r<2^{-m+1}$ and (\ref{Eq:Holder33}) is
applicable. Since $\varphi_1(r)$ is increasing near $r=0$,
(\ref{Eq:Holder3}) follows from (\ref{Eq:Holder33}).
\end{proof}

As an application of Theorem \ref{Th:Holder3}, we derive a lower bound
for the exact Hausdorff measure of the set $M_2$ of intersection times.
The corresponding problem for the upper bound remains open.

\begin{theorem} \label{Thm:dimLow}
Let $B^{\a_1}= \{B^{\a_1}(s),\,s\in\R^{N_1}\}$ and
$B^{\a_2}= \{B^{\a_2}(t),\,t\in\R^{N_2}\}$ be two independent fractional
Brownian motions with values in $\R^d$ and Hurst indices $\a_1$ and $\a_2$,
respectively.
Assume that $\frac{N_1}{\a_1}+\frac{N_2}{\a_2}>d$. Then, for
every $R > 0$, there exists a positive
constant $c_{_{4,12}}$ such that with probability 1,
\begin{equation} \label{Eq:dimLow}
\hbox{$\varphi_1$-$m$}\left(M_2\cap O_{N_1,N_2}(0,R)\right)\ge
c_{_{4,12}} L(x,O_{N_1,N_2}(0,R)),
\end{equation}
where $\varphi_1$-$m$ denotes the $\varphi_1$-Hausdorff measure.
\end{theorem}

\begin{proof}\, As in the proof of Theorem 4.1 in Xiao (1997),
(\ref{Eq:dimLow}) follows from Theorem \ref{Th:Holder3} and the
upper density theorem of Rogers and Taylor (1961). We omit the
details.
\end{proof}

As a corollary of Theorem \ref{Thm:dimLow}, we have the following
result which is stronger than (\ref{Eq:dimM20}).
\begin{corollary}\label{Cor:dimM2}
Let $B^{\a_1}$ and $B^{\a_2}$ be defined as that in Theorem \ref{Thm:dimLow}.
If $\frac{N_1}{\a_1}+\frac{N_2}{\a_2}>d$, then with probability 1,
\begin{equation}\label{Eq:dimM2}
\dim M_2 = \dimp M_2 = \left\{\begin{array}{ll}
N-\a_1d \qquad \qquad &\hbox{ if } \,  \frac{N_1} {\a_1} > d,\\
N_2+\frac{\a_2}{\a_1}N_1-\a_2d &\hbox{ if }\,  \frac{N_1}{\a_1}\le
d<\frac{N_1}{\a_1}+\frac{N_2}{\a_2}.\\
\end{array}
\right.
\end{equation}
\end{corollary}
\begin{proof}\, It is known from Theorem 7.1 in Xiao (2009) that
$\dimp M_2 \le \beta_\tau$ almost surely. In order to
prove $\dim M_2 \ge \beta_\tau$ almost surely, thanks to
Theorem \ref{Thm:dimLow}, it is sufficient to show
that with probability 1, the intersection local time $L^{\a_1, \a_2}
(O_{N_1,N_2}(0,R))>0$ for $R $ large enough. We can actually prove
a stronger result than this last statement. First note that, when $x =0$,
(\ref{Eq:m51}) becomes an equality. Thus, one can verify that
$\E[L^{\a_1, \a_2}(O_{N_1,N_2}(0,1))]>0$,
which implies that $L^{\a_1, \a_2}(O_{N_1,N_2}(0,1)) > 0$ with positive
probability. More precisely, there exist positive constants
$\delta_1$ and $\delta_2$ such that
$\P\big(L^{\a_1, \a_2}(O_{N_1,N_2}(0,1))\ge \delta_1\big)
\ge \delta_2$.

For any integer $n \ge 1$, define the event
\[
A_n = \Big\{L(0, [0, 2^{-n/\a_1}]^{N_1} \times [0, 2^{-n/\a_2}]^{N_2})
\ge \delta_1 \, 2^{-n( \frac{N_1} {\a_1} + \frac{N_2} {\a_2} - d)}\Big\}.
\]
By the scaling property (\ref{Eq:Lss}), we have $\P(A_n) \ge \delta_2$
for all $n \ge 1$.
It follows from this and Fatou's lemma that
$\P(\limsup\limits_{n\to \infty} A_n\big) \ge \delta_2$.
This implies that with positive probability
\begin{equation}\label{Eq:locallimit}
\limsup_{r \to 0} \frac{ L^{\a_1, \a_2} \big([0,\,r^{1/\a_1}]^{N_1}
\times [0,\, r^{1/\a_2}]^{N_2}) \big)}
{r^{\frac{N_1} {\a_1} + \frac{N_2}{\a_2} - d}}\ge \delta_1.
\end{equation}

Finally, note that the Gaussian field $X$ has stationary increments
and satisfies the condition of Theorem 2.1 of Pitt and Tran (1979),
which is a zero-one law for $X$ at 0. Hence (\ref{Eq:locallimit})
holds with probability 1 which, in turn,
implies $L^{\a_1, \a_2} (O_{N_1,N_2}(0,R))>0$ for all $R >0$.
\end{proof}

\section{Hausdorff and packing dimensions of $D_2$}
\label{Sec:dim}

In this section, we determine the Hausdorff and packing dimensions
of the set $D_2$ of intersection points of $B^{\a_1}$ and $B^{\a_2}$,
defined by $ D_2 = \{x \in \R^d:  x =B^{\a_1} (s)= B^{\a_2}(t) \
\hbox{ for some }\, (s, t) \in \R^N\}.$
Note that we can rewrite $D_2$ as $D_2 = B^{\a_1}(\R^{N_1}) \cap
B^{\a_2}(\R^{N_2})$.

\begin{theorem} \label{Thm:dimD2}
Let $B^{\a_1}$ and $B^{\a_2}$ be defined as that in Theorem \ref{Thm:dimLow}.
If $\frac{N_1}{\a_1}+\frac{N_2}{\a_2}>d$, then with probability 1,
\begin{equation}\label{Eq:dimD2}
\dim D_2 = \dimp D_2 = \left\{\begin{array}{ll}
d \qquad \quad &\hbox{ if } \, N_1  > \a_1 d\, \hbox{ and }\, N_2 > \a_2 d,\\
\frac{N_2} {\a_2} &\hbox{ if } \, N_1> \a_1 d\, \hbox{ and }\, N_2 \le \a_2 d,\\
\frac{N_1} {\a_1} &\hbox{ if } \, N_1 \le \a_1 d \, \hbox{ and }\, N_2 > \a_2 d,\\
\frac{N_1} {\a_1} + \frac{N_2} {\a_2} - d  &\hbox{ if } \, N_1 \le \a_1 d \,
\hbox{ and }\, N_2  \le \a_2 d.\\
\end{array}
\right.
\end{equation}
\end{theorem}

In order to prove Theorem \ref{Thm:dimD2}, we will make use of the following
two lemmas which are corollaries of the results in Monrad and Pitt (1987).

\begin{lemma}\label{Lem:Rd}
Let $B^{\a}= \{B^\a(t), t \in \R^p\}$ be a fractional Brownian motion with
values in $\R^d$ and index $\a \in (0, 1)$. If $p > \a d$, then almost surely
$B^\a(\R^p) = \R^d$.
\end{lemma}

\begin{lemma}\label{Lem:uniform}
Let $B^{\a}= \{B^\a(t), t \in \R^p\}$ be a fractional Brownian motion with
values in $\R^d$ and index $\a \in (0, 1)$. If $p \le \a d$, then for any
constants $R \ge 1$, $\ep > 0$ and $\beta > 0$ such that $0 < \alpha -
\ep < \beta <
\alpha,$ the following statement holds: With probability 1, for large
enough $n$ and for all balls $U \subseteq \R^d$ of radius $2^{-n\beta}$,
the inverse image $(B^{\a})^{-1}(U)$ can intersect at most $2^{n\ep  d}$
cubes $I_{n,\bar k}$ of the form
$$
I_{n, \bar k} = \bigl\{t \in [0,\,R]^p : (k_i - 1)2^{-n} \le t_i \le k_i2^
{-n}, \ i = 1, 2, \ldots, p \bigr\}, $$
where $\bar k = (k_1, \ldots, k_p)$ and $1 \le k_i \le R\,2^n$ for $i = 1,
\ldots, p$.
\end{lemma}

\begin{proof} {\bf of Theorem \ref{Thm:dimD2}}\ We prove (\ref{Eq:dimD2})
by considering the four cases separately.

Firstly, we assume that $N_1  > \a_1 d$ and $N_2  > \a_2 d$. It follows from Lemma
\ref{Lem:Rd} that almost surely $B^{\a_1}(\R^{N_1}) = B^{\a_2}(\R^{N_2}) = \R^d$.
Hence $D_2 = B^{\a_1}(\R^{N_1}) \cap B^{\a_2}(\R^{N_2}) = \R^d$ a.s., which
implies that $\dim D_2 = \dimp D_2 = d$ almost surely.

Secondly, we assume that $N_1  > \a_1 d$ and $N_2  \le \a_2 d$. Then $D_2 =
B^{\a_1}(\R^{N_1}) \cap B^{\a_2}(\R^{N_2}) = B^{\a_2}(\R^{N_2})$ a.s., which yields
$\dim D_2 = \dimp D_2 = N_2/\a_2$ almost surely. The proof for the case
$N_1  \le \a_1 d$ and $N_2 > \a_2 d$ is similar.

Finally, we consider the case of $N_1 \le \a_1 d$ and $N_2  \le \a_2 d$
[in addition to $\frac{N_1}{\a_1}+\frac{N_2}{\a_2}>d$].
Let $S_2$ be the projection of $M_2$ on $\R^{N_2}$. Then $B^{\a_2}(S_2) = D_2.$
Since, for every $\ep > 0$, $B^{\a_2}(t)$ satisfies a uniform H\"older
condition of order
$\alpha_2 - \ep$ on every compact interval of $\R^{N_2}$, we have
\begin{equation}\label{Eq:dimD2-3}
\dimp B^{\a_2}(S_2) \le \frac 1 {\alpha_2}\, \dimp S_2 \le \frac{N_1}
 {\a_1} + \frac{N_2} {\a_2} - d ,\quad
a.s.,
\end{equation}
where the last inequality follows from (\ref{Eq:dimM2}).

It only remains to prove $\dim D_2 \ge \frac{N_1} {\a_1} + \frac{N_2}
{\a_2} - d$ almost surely. For this purpose, we denote
$\ell = \frac{N_1} {\a_1} + \frac{N_2} {\a_2} - d$ and define an
$(N, 2d)$-Gaussian random field $Z = \{Z(s, t), (s, t) \in \R^N\}$ by
$$
Z(s, t) = \big(B^{\a_1}(s), B^{\a_2}(t) \big),\quad  (s, t) \in \R^N.
$$
Set $ \tilde D_2 = \{(x, x): x \in D_2\}$. Then
\begin{equation}\label{Eq:dimD2-4}
Z^{-1}(\tilde D_2) = M_2.
\end{equation}

Fix an $\omega \in \Omega$ such that the conclusion of Lemma
\ref{Lem:uniform} holds. Assume that for some constant $\eta > 0$,
$\dim D_2(\omega) < \ell-\eta$. [We will suppress $\omega$ from now
on.] Then, for any $n$ large enough, there exists a sequence of
balls $\{U_i\}$ in $\R^d$ with radius $\le 2^{-n},$ such that
\begin{equation}\label{Eq:dimD2-5}
D_2 \subseteq \bigcup_i
U_i \quad \hbox{ and }\ \sum_i ({\rm diam} U_i)^{\ell-\eta} \le 1,
\end{equation}
where diam$U$ denotes the diameter of $U$. Choose positive constants
$\ep$, $\ga_1 < \a_1$ and $\gamma_2< \a_2$ such that
\begin{equation}\label{Eq:const-ep}
\ga_1 < \ga_2 \quad \hbox{ and }\ \ep d \Big(\frac{\ga_2} {\ga_1} +
1\Big) < \frac {\ga_2 \eta}  2.
\end{equation}

Let $m_i$ and $ n_i$  be integers that satisfy
\begin{equation}\label{Eq:dimD2-6}
2^{-(m_i+1)\ga_1} \le {\rm diam} U_i \le 2^{-m_i\ga_1} \ \ \ \hbox{
and }\ \ 2^{-(n_i+1)\ga_2} \le {\rm diam} U_i \le 2^{-n_i\ga_2}.
\end{equation}
By (\ref{Eq:dimD2-4}) and Lemma \ref{Lem:uniform}, we have
\begin{equation}\label{Eq:dimD2-7}
\begin{split}
M_2 \cap [0, R]^N &\subseteq \bigcup_i Z^{-1}(U_i\times U_i) \\
&\subseteq \bigcup_i \Big\{\text{the union of at most
$2^{(m_i+ n_i)\ep d}$ cubes } I_{m_i, \bar j} \times I_{n_i, \bar k}\Big\}.
\end{split}
\end{equation}
Denote the cubes in the right-hand side of (\ref{Eq:dimD2-7}) by
$C_{ij}.$ Note that, since $\ga_1 < \ga_2$, we derive from
(\ref{Eq:dimD2-6}) that diam$C_{ij}\le 3 \cdot 2^{-n_i}$ for $i$ (or
$n$) large enough. Combining this with (\ref{Eq:dimD2-5}),
(\ref{Eq:dimD2-6}) and (\ref{Eq:dimD2-6}), we derive that for all
$n$ large enough,
\begin{equation}\label{Eq:dimD2-8}
\begin{split}
\sum_i\sum_j ({\rm diam} C_{ij})^{\ga_2 (\ell- \frac {\eta} 2)}
&\le c_{_{5, 1}}\, \sum_i 2^{(m_i+ n_i)\ep d}(2^{-n_i})^{\ga_2 (\ell-\frac {\eta} 2)}\\
&\le c_{_{5, 2}}\, \sum_i 2^{n_i(1+ \frac{\ga_2} {\ga_1})\ep d } (2^{-n_i \ga_2})^{\ell- \frac \eta 2}\\
&\le c_{_{5, 3}}\,\sum_i ({\rm diam} U_i)^{\ell-\eta} \le c_{_{5, 3}}.
\end{split}
\end{equation}
It follows from (\ref{Eq:dimD2-7}) and  (\ref{Eq:dimD2-8}) that
$$
\dim \big(M_2\cap [0, R]^N\big) \le \ga_2 \bigl(\ell - \frac {\eta} 2 \bigr) <
\alpha_2\bigl(\ell - \frac {\eta} 2 \bigr).
$$
Hence we have proven that, for any $\eta > 0$,
\begin{equation}\label{Eq:dimD2-9}
\P\Big\{ \dim D_2 \ge \ell-\eta \Big\} \ge \P\Big\{\dim \big(M_2 \cap [0, R]^N\big)\ge \alpha_2
\bigl(\ell - \frac {\eta} 2 \bigr)\Big\}.
\end{equation}
Letting $R \uparrow \infty$ and $\eta \downarrow 0$ along the rational numbers
and by using (\ref{Eq:dimM2}),
we obtain $\dim D_2 \ge \ell$ almost surely. This finishes the proof.
\end{proof}

\bibliographystyle{plain}
\begin{small}

\end{small}

\vspace{.2in}

\begin{quote}
\begin{small}

 \noindent \textsc{Dongsheng Wu}.\
        Department of Mathematical Sciences, 201J Shelby Center,
        University of Alabama in Huntsville,
        Huntsville, AL 35899, U.S.A.\\
        E-mail: \texttt{dongsheng.wu@uah.edu}\\
        URL: \texttt{http://webpages.uah.edu/\~{}dw0001}\\

  \noindent \textsc{Yimin Xiao}.\
        Department of Statistics and Probability, A-413 Wells
        Hall, Michigan State University,
        East Lansing, MI 48824, U.S.A.\\
        E-mail: \texttt{ xiao@stt.msu.edu}\\
        URL: \texttt{http://www.stt.msu.edu/\~{}xiaoyimi}
\end{small}
\end{quote}
\end{document}